\documentclass{article}
\usepackage{amssymb}
\usepackage{amsmath}
\usepackage{algorithm}
\usepackage{algorithmic}
\usepackage{amsthm}
\usepackage{booktabs}
\usepackage{threeparttable}
\usepackage{graphicx}
\usepackage{cite}
\usepackage{listings}
\usepackage{xcolor}
\usepackage{epstopdf}
\usepackage{multirow}
\newtheorem{Lemma}{Lemma}[section]
\newtheorem{Theorem}{Theorem}[section]
\newtheorem{Definition}{Definition}[section]
\newtheorem{Proposition}{Proposition}[section]

\begin{document}

\title{An adaptive cubic regularization method for computing extreme eigenvalues of tensors}
\author{
  {Jingya Chang}\thanks{%
    School of  Mathematics and Statistics, Guangdong University of Technology
    ({\tt jychang@gdut.edu.cn}). This author's work was supported by the National Natural Science Foundation of China (grant No. 11901118 and 62073087).}
    \ and
     {Zhi Zhu}\thanks{%
    School of Mathematics and Statistics, Guangdong University of Technology    ({\tt zhuzhixixi@163.com}).}
  }

  \date{}
  \maketitle

\noindent{\textbf{Abstract.}} In this paper, we compute the H- and Z-eigenvalues of even order symmetric tensors by using the adaptive cubic regularization algorithm. First, the equation  of eigenvalues of  the tensor  is represented by  a spherically constrained optimization problem. Owing to the nice geometry of the spherical constraint, we minimize the objective function and preserve the constraint in an alternating way.
  The objective function of our optimization model is approximated by a cubic function with a tunable parameter, which is solved  inexactly to obtain a trial step. Then the Cayley transform is applied to the trial step. Based on the ratio of actual and predicted reductions, a parameter  is regulated to make sure that the cubic function is a good estimation of the original objective function. Finally we obtain our adaptive cubic regularization algorithm for computing an eigenvalue of a tensor (ACRCET). Furthermore, we prove that the sequence of iterations  generated by  ACRCET  converges to an eigenvalue of a given tensor globally. In order to improve the computational efficiency, we propose a fast computing skill for $\mathcal{T}\mathbf{x}^{r-2}$ which is the matrix-valued product of a hypergraph related tensor $\mathcal{T}$ and a vector $\mathbf{x}$. Numerical experiments illustrate that the fast computing skill for $\mathcal{T}\mathbf{x}^{r-2}$ is efficient and our ACRCET is effective  when computing eigenvalues of   even order symmetric tensors. \\




\section{Introduction}
Eigenvalues of tensors were proposed by Qi \cite{qi2005eigenvalues} and Lim \cite{lim2005singular} respectively in 2005. From then on, eigenvalues and eigenvectors of tensors have been widely  used in science and engineering, such as medical imaging, image processing and spectral graph theory. For example, in magnetic resonance imaging, the Z-eigenvectors of the major Z-eigenvalues of an even order tensor express the orientations of crossing nerve fibers in white matter of human brain \cite{chen2013positive}. The limiting probability distribution vector of a higher-order Markov chain is a Z-eigenvector of a transition probability tensor, and the corresponding eigenvalue is 1  \cite{li2014limiting}. The Z-eigenvalues of adjacency tensors of even-uniform hypergraphs also have important applications in spectral hypergraph theory \cite{li2013z}.

In recent years, the computation of eigenvalues of tensors has been studied by many scholars and different kinds of methods have been proposed.
For nonnegative tensors, Ng and Qi \cite{ng2010finding} extended the  method of Collatz (1942) for calculating the spectral radius of an irreducible nonnegative matrix to the calculation of the largest eigenvalue of a nonnegative tensor. A homotopy method was proposed to compute the largest eigenvalues of irreducible nonnegative tensors in \cite{chen2019homotopy}. Kuo  et al. gave a   homotopy continuation method for the computation of nonnegative Z-/H-eigenpairs of nonnegative tensors \cite{kuo2018continuation}.
 Kolda and Mayo \cite{kolda2011shifted} proposed a shifted power method to compute Z-eigenpairs of symmetric tensors. Later,  they developed  an adaptive shifted power method for computing generalized tensor Z-eigenpairs \cite{kolda2014adaptive}. A shifted inverse power method was introduce in \cite{sheng2021computing}  for computing Z-eigenvalues of tensors.  Han \cite{han2013unconstrained}  introduced two unconstrained minimization models to compute the Z-, H-, or D-eigenvalues of even order symmetric tensors.

For computing all eigenvalues of a tensor, Cui et al. \cite{cui2014all} used the Jacobian semidefinite relaxations in polynomial optimization to calculate the eigenvalues of a tensor sequentially from the largest one to the smallest one. Chen et al. \cite{chen2016computing} gave two  homotopy continuation type algorithms which can find all equivalent classes of isolated generalized eigenpairs if executed properly.
For large scale tensors, Chen et al. \cite{chenq2016computing} introduced an algorithm for computing extreme eigenvalues of large scale Hankel tensors by using the fast Fourier transform.
Chang et al. \cite{chang2016computing} proposed a limited memory BFGS quasi-Newton algorithm for computing eigenvalues of large scale sparse tensors related with a hypergraph.

The cubic estimation model was first introduced in \cite{griewank1981modification} for finding an improved Newton step. Later, the convergence property and numerical performance of the cubic regularization method were broadly studied \cite{cartis2011adaptive,nesterov2006cubic,xu2020newton}.  Specifically, some smart worst-case global iteration complexity bounds were established in \cite{cartis2011adaptive2}.
 The cubic regularization method provides another option beyond trust region and line search methods for unconstrained optimization problems. In this paper, we generalize the cubic regularization  method from unconstrained optimization to spherically constrained optimization.

For even order symmetric tensors, we convert the tensor eigenvalue problems to a spherically constrained minimization problem equivalently. To solve this orthogonally constrained model, first we  surrogate the objective function by a third order model with  an approximated Hessian and find the minimal point of the third order approximation  imprecisely at each iteration. Next we push the minimal point onto the the unit sphere by Cayley transform to keep the iteration point feasible. In this step, we use a curvilinear line search so that the objective function value has an appropriate decrease. The approximation function  contains a parameter which is adjusted during the iterative process to insure a good estimation. In terms of eigenvalue problem of tensors arising from a hypergraph, we  propose a fast computation algorithm  for the repeatedly occurred operation $\mathcal{T}\mathbf{x}^{r-2},$ which is the product of a vector and a tensor arising from a hypergraph. Moreover, we analyze the convergence property of the iterative sequence and prove that the iteration points converge to an eigenvalue of the tensor. We perform numerical experiments and compare our method with PM \cite{kolda2011shifted,kolda2014adaptive}, ACSA \cite{chenq2016computing} and HUOA \cite{han2013unconstrained} for computing eigenvalues of even order symmetric tensors.

The outline of this paper is drawn as follows. Basic knowledge about  tensors and tensor eigenvalues is given in Section 2. In Section 3, we introduce our algorithm for computing eigenvalues and eigenvectors of an even order symmetric tensor. In Section 4, we demonstrate the global convergence property of our  method. In Section 5, we show the fast computation technique for  $\mathcal{T}\mathbf{x}^{r-2}$. Numerical experiments on small and medium scale tensors are given in Section 6. Finally, we give some concluding remarks in Section 7.

\section{Preliminary}
 We use boldface Euler script letters such as $\mathcal{A}$ to represent  tensors. A matrix is named with a capital letter, while a lower case bold letter is used for a vector and a lower case letter for a scalar.  Denote $\mathbb{R}^{[r,n]}$ as the space of $r$th order $n$-dimensional real tensors, $\mathbb{R}^{m\times n}$ as the space of  real matrices with $m$ rows and $n$ columns, and $\mathbb{R}^n$ as the space of $n$-dimensional real vectors, where $r$, $m$ and $n$ are positive integers. Before going to the main results, we introduce the concepts related to tensors.

 A tensor $\mathcal{A}\in\mathbb{R}^{[r,n]}$ has $n^r$ entries:
\[\{a_{i_1i_2\cdots i_r}\} \] \text{for} $  i_j\in \{1,2,\cdots,n\} \ \ \text{and} \ \ j\in \{1,2,\cdots,r\}.$
A tensor $\mathcal{A}\in\mathbb{R}^{[r,n]}$ is a symmetric tensor if the value of $a_{i_1i_2\cdots i_r}$ is invariable under any permutation of its indices. An identity tensor $\mathcal{I}\in\mathbb{R}^{[r,n]}$ is a tensor   whose diagonal entries are all one and other off-diagonal entries are zero.

For a vector $\mathbf{x}\in\mathbb{R}^n$, we define a scalar $\mathcal{A}\mathbf{x}^r\in\mathbb{R}$,
\[
\mathcal{A}\mathbf{x}^r=\sum_{i_1,\cdots,i_r=1}^{n}a_{i_1i_2i_3\cdots i_r}x_{i_1}x_{i_2}x_{i_3}\cdots x_{i_r},
\]
a vector $\mathcal{A}\mathbf{x}^{r-1}\in\mathbb{R}^n$,
\[
(\mathcal{A}\mathbf{x}^{r-1})_p=\sum_{i_2,\cdots,i_r=1}^{n}a_{pi_2i_3\cdots i_r}x_{i_2}x_{i_3}\cdots x_{i_r}, \forall p\in\{1,2,\cdots,n\},
\]
and a matrix $\mathcal{A}\mathbf{x}^{r-2}\in\mathbb{R}^{n\times n}$,
\[
(\mathcal{A}\mathbf{x}^{r-2})_{pq}=\sum_{i_3,\cdots,i_r=1}^{n}a_{pqi_3\cdots i_r}x_{i_3}\cdots x_{i_r}, \forall p,q\in\{1,2,\cdots,n\}.
\]
\begin{Definition}[H-eigenvalue and H-eigenvector \cite{qi2005eigenvalues}]
If there exist a $\lambda\in\mathbb{R}$ and a nonzero vector $\mathbf{x}\in\mathbb{R}^n$ such that
\begin{equation}\label{H-eigenvalue}
\mathcal{A}\mathbf{x}^{r-1}=\lambda\mathbf{x}^{[r-1]},
\end{equation}
then $\lambda$ is called an H-eigenvalue of $\mathcal{A}$ and $\mathbf{x}$ is its associated H-eigenvector, where $\mathbf{x}^{[r-1]}=:(x^{r-1}_1,x^{r-2}_2,\cdots,x^{r-1}_n)^T$.
\end{Definition}
\begin{Definition}[Z-eigenvalue and Z-eigenvector \cite{qi2005eigenvalues}]
Suppose that $(\lambda,\mathbf{x})$ is a solution pair of the following system
\begin{equation}\label{Z-eigenvalue}
\mathcal{A}\mathbf{x}^{r-1}=\lambda\mathbf{x}\qquad {\rm and} \qquad\mathbf{x}^T\mathbf{x}=1,
\end{equation}
then $\lambda$ is called a Z-eigenvalue of $\mathcal{A}$ and $\mathbf{x}$ is its associated Z-eigenvector.
\end{Definition}

\section{Computation of eigenvalues of tensors}
In this section, we demonstrate our method for computing H- and Z-eigenvalues of an even order symmetric tensor.
\subsection{The equivalent optimization model of tensor eigenvalue problems}
The equation systems in \eqref{H-eigenvalue} and \eqref{Z-eigenvalue} can be changed to a  spherically constrained optimization problem \cite{chang2016computing}. Consider the following optimization problem
\begin{equation}\label{OptModel}
\min \ f(\mathbf{x}):=\frac{\mathcal{A}\mathbf{x}^r}{\mathcal{B}\mathbf{x}^r} \qquad {\rm s.t.} \ \  \mathbf{x}\in\mathbb{S}^{n-1},
\end{equation}
where  $\mathcal{A}$ and  $\mathcal{B}$ are symmetric tensors, the set $\mathbb{S}^{n-1}:=\{\mathbf{x}\in \mathbb{R}^n:\mathbf{x}^T\mathbf{x}=1\}$ is a unit spherical surface. The gradient of $f(\mathbf{x})$ is
\begin{equation} \label{gradient}
\nabla f(\mathbf{x})=\frac{r}{\mathcal{B}\mathbf{x}^r}\left(\mathcal{A}\mathbf{x}^{r-1}-\frac {\mathcal{A}\mathbf{x}^r}{\mathcal{B}\mathbf{x}^r}\mathcal{B}\mathbf{x}^{r-1}\right
).
\end{equation}
When $\nabla f(\mathbf{x})=0,$ we have $\mathcal{A}\mathbf{x}^{r-1}=\frac {\mathcal{A}\mathbf{x}^r}{\mathcal{B}\mathbf{x}^r}\mathcal{B}\mathbf{x}^{r-1}.$
In \cite{han2013unconstrained}, it is shown that when $\mathcal{B}=\mathcal{I},$
$$\mathcal{B}x^{r}=\|\mathbf{x}\|_r^r, \ \ \mathcal{B}x^{r-1}=\mathbf{x}^{[r-1]}.$$
When $r$ is even and $\mathcal{B}=I^{\frac{r}{2}},$
$$\mathcal{B}\mathbf{x}^{r}=\|\mathbf{x}\|_2^r=1, \ \ \mathcal{B}\mathbf{x}^{r-1}=\| \mathbf{x}\|_2^{r-2} \mathbf{x}= \mathbf{x}.$$
Thus for even order tensors,
$f(\mathbf{x})=\mathcal{A}\mathbf{x}^r/\mathcal{B}\mathbf{x}^r$ is the  H-eigenvalue or Z-eigenvalue of $\mathcal{A}$ when   $\mathcal{B}=\mathcal{I}$ or $\mathcal{B}=I^{\frac{r}{2}}$ respectively at the stationary point $\mathbf{x},$ and $\mathbf{x}$ is the corresponding eigenvector. On the other hand, we can verify that when  $\mathcal{B}=\mathcal{I}$ or $ \mathcal{B}=I^{\frac{r}{2}}$ and $\mathbf{x}$ is the H- or Z-eigenvector of $\mathcal{A},$ the eigenvector $\mathbf{x}$ is the stationary point with the eigenvalue being $\frac {\mathcal{A}\mathbf{x}^r}{\mathcal{B}\mathbf{x}^r}.$

Therefore, the eigenvalue problems are equivalently transformed to the question of finding the stationary point of the objective function in \eqref{OptModel} on the unit sphere.

\subsection{Adaptive cubic regularization method}
In this section, we  design an adaptive cubic regularization method for solving the spherically constrained optimization problem \eqref{OptModel}.
 By calculating,  we get the Hessian $\nabla^2f(\mathbf{x})$ of $f(\mathbf{x})$ as follows
 \begin{align}\label{Hessian}
\nabla^2f(\mathbf{x}_k) =&\frac{r(r-1)}{\mathcal{B}\mathbf{x}^r}\mathcal{A}\mathbf{x}^{r-2}-\frac{r^2}{(\mathcal{B}\mathbf{x}^r)^2}(\mathcal{A}\mathbf{x}^{r-1}\circledcirc \mathcal{B}\mathbf{x}^{r-1})\\
                             &- \frac{r(r-1)\mathcal{A}\mathbf{x}^r}{(\mathcal{B}\mathbf{x}^r)^2}\mathcal{B}\mathbf{x}^{r-2}+\frac{r^2\mathcal{A}\mathbf{x}^r}{(\mathcal{B}\mathbf{x}^r)^3}(\mathcal{B}\mathbf{x}^{r-1}\circledcirc \mathcal{B}\mathbf{x}^{r-1}).
\end{align}
The symbol $\circledcirc$ refers to the operation $\mathbf{a}\circledcirc\mathbf{b}=\mathbf{ab}^T+\mathbf{ba}^T$ for two vectors $\mathbf{a}$ and $\mathbf{b}$. Clearly, the function $f(\mathbf{x})$ is three times continuously differentiable.
\begin{Lemma}\label{Upperbound}
Since $f(\mathbf{x})$ is three times continuously differentiable on the the compact  set $\mathbb{S}^{n-1},$
there exists a constant $M$ such that
\begin{equation}\label{Bound}
 \|f(\mathbf{x})\|\leq M, \|\nabla f(\mathbf{x})\|\leq M, \|\nabla^2f(\mathbf{x})\|\leq M.
\end{equation}
Also $\nabla^2f(\mathbf{x})$ is globally Lipschitz continuous, which means  there exists a constant $L>0$ such that
\begin{equation}\label{Lipschitz}
  \| \nabla^2f(\mathbf{x})-\nabla^2f(\mathbf{\tilde x})\| \leq L\|\mathbf{x}-\mathbf{\tilde x}\|
\end{equation}
 for any $\mathbf{x}$ and $\mathbf{\tilde x} \in \mathbb{S}^{n-1}$.
\end{Lemma}

Denote $f_k=f(\mathbf{x}_k)$, $\mathbf{g}_k=\nabla f(\mathbf{x}_k)$ and $H_k=\nabla^2f(\mathbf{x}_k).$  Let $\mathbf{p}$ be a direction
pointing from the vector $\mathbf{x}_k.$
By expressing the Taylor expansion of $f(\mathbf{x}_k+\mathbf{p})$ at the point $\mathbf{x}_k,$ we get
\begin{align}
f(\mathbf{x}_k+\mathbf{p}) &= f_k+\mathbf{g}_k^T\mathbf{p}+\frac{1}{2}\mathbf{p}^TH_k\mathbf{p}+\int_0^1(1-t)\mathbf{p}^T[H(\mathbf{x}_k+t \mathbf{p})-H(\mathbf{x}_k)]\mathbf{p}\ \mathrm{d}t \notag \\
                           &\leq f_k+\mathbf{g}_k^T\mathbf{p}+\frac{1}{2}\mathbf{p}^TH_k\mathbf{p}+\frac{1}{6}L\|\mathbf{p}\|_2^3,
\end{align}
in which the inequality is deduced from  the Lipschitz property of $\nabla^2 f(\mathbf{x}).$
In \cite{cartis2011adaptive}, Cartis et al. introduced  a dynamic parameter $\sigma_k $ instead of $\frac{1}{2}L$ and suggested to use
\begin{equation}\label{AppMod0}
 f_k+\mathbf{g}_k^T\mathbf{p}+\frac{1}{2}\mathbf{p}^TH_k\mathbf{p}+\frac{1}{3}\sigma_k\|\mathbf{p}\|_2^3
\end{equation}
as the approximation to $f(\mathbf{x}_k+\mathbf{p}).$ Since the function $f$ is restricted on the unit spherical surface, we project the Hessian matrix $H_k$ onto the tangent space at $\mathbf{x}$ by the projection matrix $P_k:=(I-\mathbf{x}_k\mathbf{x}_k^T),$ and get
$B_k:=P_kH_kP_k.$ The estimation in \eqref{AppMod0} is further modified by
\begin{equation}
 f_k+\mathbf{g}_k^T\mathbf{p}+\frac{1}{2}\mathbf{p}^TB_k\mathbf{p}+\frac{1}{3}\sigma_k \|\mathbf{p}\|_2^3.
\end{equation}
Thus in each iteration, we solve the subproblem
\begin{equation}
\min \ m_k(\mathbf{p}):=f_k+\mathbf{g}_k^T\mathbf{p}+\frac{1}{2}\mathbf{p}^TB_k\mathbf{p}+\frac{1}{3}\sigma_k\|\mathbf{p}\|_2^3 \label{sub}
\end{equation}
 to find a descent direction $\mathbf{p}_k$.

 In order to improve the computing efficiency, the  subproblem \eqref{sub} is solved inexactly.  Denote $\mathbf{p}_k^C$ as the Cauchy point of $ m_k(\mathbf{p}),$ that means
\begin{equation*}
  \mathbf{p}_k^C = -\tau_k^C \mathbf{g}_k \quad \text{and} \quad \tau_k^C= \arg\min_{\tau\in R^+} m_k(-\tau\mathbf{g}_k).
\end{equation*}A  vector $\mathbf{p}_k$ is chosen once it satisfies
 \begin{equation*}
   m_k(\mathbf{p}_k)\leq m_k(\mathbf{p}_k^C).
 \end{equation*}
 The parameter $\sigma_k$ is updated adaptively in  light of the ratio of actual reduction to predicted reduction, details of which will be explained in the next subsection.
\subsection{Cayley transform, curvilinear  search and parameter tuning  }
 Based on the direction $\mathbf{p}_k$ getting from \eqref{sub}, we employ the Cayley transform \cite{cayley1846algebraic} to generate an orthogonal matrix and obtain a new feasible point. At the same time, we perform a curvilinear search on  $\mathbb{S}^{n-1}$  to guarantee  that the objective function value decrease. The parameter $\sigma$ is tuned according to the accuracy of the estimation term.

Define a skew-symmetric matrix $W_k(\alpha)\in \mathbb{R}^{n \times n}$ as
\begin{equation}\label{SkewSymMatr}
  W_k(\alpha)=\frac{\alpha}{2}(\mathbf{x}_k\mathbf{p}^T_k-\mathbf{p}_k\mathbf{x}^T_k),
\end{equation}
where $\alpha>0$ is a parameter. Then $(I+W_k(\alpha))$ is invertible, and the Cayley transform produces an orthogonal matrix
\begin{equation}\label{OrthogMatr}
  Q_k(\alpha)=(I+W_k(\alpha))^{-1}(I-W_k(\alpha)).
\end{equation}
For $\mathbf{x}_k \in \mathbb{S}^{n-1},$ the vector
\begin{equation}\label{xk+0}
  \mathbf{x}_k^+(\alpha):=Q_k(\alpha)\mathbf{x}_k
\end{equation}
also belongs to $\mathbb{S}^{n-1}.$
Substituting \eqref{SkewSymMatr} and \eqref{OrthogMatr} into \eqref{xk+0}, the vector $\mathbf{x}_k^+(\alpha)$ can be computed through the following formula \cite{chang2016computing}
\begin{equation}
\mathbf{x}_k^+(\alpha)=\frac{[(2-\alpha \mathbf{p}^T_k \mathbf{x}_k)^2-\alpha^2\|\mathbf{p}_k\|^2]\mathbf{x}_k+4\alpha \mathbf{p}_k}{4+\alpha^2\|\mathbf{p}_k\|^2-\alpha^2(\mathbf{p}_k^T\mathbf{x}_k)^2}.\label{xk+}
\end{equation}

Next, a  curvilinear search process by \eqref{xk+}  is implemented to determine the step
size $\alpha$.  Define the ratio of actual reduction in  $f$ to  predicted reduction in  $m_k$ as
\begin{equation}
\rho_k:=\frac{f(\mathbf{x}_k)-f(\mathbf{x}_k^+(\alpha))}{m_k(0)-m_k(\alpha\mathbf{p}_k)}.
\end{equation}
Take $\alpha = 1$ initially. If  $\rho_k$ is greater than or equal to a positive constant $\eta_1$($\eta_1 \in (0,1)$), the step size $\alpha$ is accepted and the new iterate $\mathbf{x}_{k+1}=\mathbf{x}_k^+(\alpha)$. Otherwise, we decrease the step size $\alpha$ and repeat until $\rho_k \geq \eta_1.$
If $\alpha_k=1$ and $\rho_k$ is greater than or equal to a given number, we regard $m_k$  as a very successful estimation of the original function $f$ and decrease the value of $\sigma$ in the next step. If not,  the value of $\sigma$ is increased in the next step. We describe the updating scheme  in detail in Algorithm \ref{algorithm1}.

\begin{algorithm}
\caption{An adaptive cubic regularization algorithm for computing an eigenvalue of a tensor.(ACRCET)}
\label{algorithm1}
\begin{algorithmic}[1]
\STATE Set $\mathcal{B}=\mathcal{I}$ and $\mathcal{B}=I^\frac{r}{2}$ when calculating an H-eigenvalue and a Z-eigenvalue of a tensor $\mathcal{A}$ respectively.
\STATE Set parameters $\gamma_3\geq \gamma_2>1>\gamma_1>0$, $1>\eta_2\geq \eta_1>0$, and $\sigma_0>0$ for $k=0,1,\cdots$ until convergence. Choose an initial point $\mathbf{x}_0\in \mathbb{S}^{n-1}$ and set $k\gets0$.
\WHILE{$\nabla f(\mathbf{x}_k)\neq 0$}
\STATE Calculate $\mathcal{A}\mathbf{x}^r, \mathcal{B}\mathbf{x}^r,\mathcal{A}\mathbf{x}^{r-1},\mathcal{B}\mathbf{x}^{r-1},\mathcal{A}\mathbf{x}^{r-2}, \text{and} \  \mathcal{B}\mathbf{x}^{r-2}.$
\STATE Solve the subproblem:
 \[\min \ f_k+\mathbf{g}_k^T\mathbf{p}+\frac{1}{2}\mathbf{p}^TB_k\mathbf{p}+\frac{1}{3}\sigma_k\|\mathbf{p}\|^3,\]
for a trial step $\mathbf{p}_k$ satisfying
\[m_k(\mathbf{p}_k )\leq m_k(\mathbf{p}_k^C)\]
where the Cauchy point
\[\mathbf{p}_k^C=-\tau_k^C\mathbf{g}_k \ \text{and} \ \tau_k^C=\arg \min_{\tau \in \mathbb{R}^+} m_k(-\tau \mathbf{g}_k).\]
\STATE Compute $f(\mathbf{x}_k^+(\alpha))$ by \eqref{xk+} and \eqref{OptModel}. Find the smallest nonnegative integer $j$  such that the step size $\alpha=\gamma_1^j$ satisfies
\[\rho_k=\frac{f(\mathbf{x}_k)-f(\mathbf{x}_k^+(\alpha))}{m_k(0)-m_k(\alpha \mathbf{p}_k)}\geq \eta_1.\]
\STATE Set $\alpha_k=\gamma_1^j$ and $\mathbf{x}_{k+1}=\mathbf{x}_k^+(\alpha_k)$.
\STATE If $\alpha_k=1$, set
\[\sigma_{k+1} \in \begin{cases}
\left[0,\sigma_k \right]                       & \text{if }\rho_k>\eta_2   \  \qquad \qquad \text{[very successful iteration]},\\
\left[\sigma_k,\gamma_2\sigma_k \right]        & \text{if }\eta_1\leq\rho_k\leq\eta_2
\qquad \text{[successful iteration]},
\end{cases}\]
\text{else}
\[\sigma_{k+1} \in [\gamma_2\sigma_k,\gamma_3\sigma_k]  \qquad \qquad \qquad  \qquad \text{[unsuccessful iteration]}.
\]

\STATE Set $k\gets{k+1}$

\ENDWHILE

\end{algorithmic}
\end{algorithm}

\section{Convergence analysis}
In this section, we analyze the convergence property of the ACRCET method. We show that the gradient norm $\|\mathbf{g}_k\|$ generated by Algorithm \ref{algorithm1} converges to $0$ globally.  Thus we obtain the stationary point which is an eigenvector of the given tensor, and its function value is the  eigenvalue.

First we explain that the step size $\alpha_k$ is bounded away from zero  in our iteration process.   Then we prove
 \begin{equation}\label{liminfg_k}
  \liminf_{k\rightarrow \infty}\|\mathbf{g}_k\|=0
 \end{equation}
by contradiction.
Based on \eqref{liminfg_k}, we further obtain the convergence result of $\{\mathbf{g}_k\}$.

Cartis et al. \cite{cartis2011adaptive} gave a  lower bound on the decrease in $f$ predicted from the cubic model $m_k(\mathbf{p}_k)$ and a bound on  $\Vert \mathbf{p}_k \Vert$. These two conclusions, which are useful in the convergence analysis, also hold for our algorithm.
\begin{Lemma}[\cite{cartis2011adaptive}] \label{CartisLm}
Suppose that the step $\mathbf{p}_k$ satisfies $m_k(\mathbf{p}_k)\leq m_k(\mathbf{p}_k^C)$. Then for all $k \geq 0$, we have
\begin{equation}\label{CartisLm1}
f(\mathbf{x}_k)-m_k(\mathbf{p}_k) \geq \frac{\Vert \mathbf{g}_k \Vert}{6\sqrt{2}}\min\ \left[\frac{\Vert \mathbf{g}_k \Vert}{1+\Vert B_k \Vert},\frac{1}{2} \sqrt{\frac{\Vert \mathbf{g}_k\Vert}{\sigma_k}}\right]
\end{equation}
and
\begin{equation}\label{CartisLm2}
\Vert \mathbf{p}_k \Vert \leq \frac{3}{\sigma_k}\max\ (M,\sqrt{\sigma_k \Vert \mathbf{g}_k \Vert} ).
\end{equation}
\end{Lemma}
Now we prove that $\alpha_k$ is bounded above $0$.
\begin{Lemma}\label{alphakBound}
For the step size $\alpha_k$ generated by Algorithm \ref{algorithm1}, it holds that
\[
\liminf_{k\to\infty}\alpha_k>0.
\]
\end{Lemma}
\begin{proof}
We prove the conclusion by contradiction. Suppose that a subsequence of   $\{\alpha_{k_i}\}$ tends to 0. By the rule of backtracking search in Algorithm \ref{algorithm1}, we have
\begin{equation}\label{L4.50}
 \frac{f_{k_i}-f(x^+_{k_i}(\gamma_1^{-1}\alpha_{k_i}))}{f_{k_i}-m_{k_i}(\gamma_1^{-1}\alpha_{k_i}\mathbf{p}_{k_i})}<\eta_1.
\end{equation}
Since the objective function is zero-order homogeneous, then
\begin{equation}\label{x_kg_k=0}
  \mathbf{g}_k^T\mathbf{x}_k=0.
\end{equation}
It can be calculated from  \eqref{xk+} that
\begin{equation}\label{x+'}
  \mathbf{x}^{'+}_k(0)=-\mathbf{x}_k^T\mathbf{p}_k\mathbf{x}_k+\mathbf{p}_k.
\end{equation}
By \eqref{x_kg_k=0} and \eqref{x+'} we obtain
\[\left.
\frac{\mathrm{d}f(\mathbf{x}_k^+(\alpha))}{\mathrm{d}\alpha}\right|_{\alpha=0}=\nabla f(\mathbf{x}_k^+(0))^T\mathbf{x}^{'+}_k(0)=\nabla f(\mathbf{x}_k)^T(-\mathbf{x}_k^T\mathbf{p}_k\mathbf{x}_k+\mathbf{p}_k)=\mathbf{g}_k^T\mathbf{p}_k
\]
and $$\frac{\mathrm{d}^2f(\mathbf{x}_k^+(\alpha))}{\mathrm{d}\alpha^2}|_{\alpha=0}=\mathbf{p}_k^TB_k\mathbf{p}_k.$$
Thus by substituting the  Taylor series expansion of $f(x^+_{k_i}(\gamma_1^{-1}\alpha_{k_i}))$ around the point $\alpha=0$ and the expression of $m_{k_i}(\gamma_1^{-1}\alpha_{k_i}\mathbf{p}_{k_i})$ into \eqref{L4.50}, we have

\begin{align*}
\eta_1(-\gamma_1^{-1}\alpha_{k_i}\mathbf{g}_{k_i}^T\mathbf{p}_{k_i}
-\frac{1}{2}(\gamma_1^{-1}\alpha_{k_i})^2\mathbf{p}_{k_i}^TB_{k_i}\mathbf{p}_{k_i}
-\frac{1}{3}\sigma_k\Vert \gamma_1^{-1}\alpha_{k_i}\mathbf{p}_{k_i}\Vert^3)&\\
+\gamma_1^{-1}\alpha_{k_i}\mathbf{g}_{k_i}^T\mathbf{p}_{k_i}
+\frac{1}{2}(\gamma_1^{-1}\alpha_{k_i})^2\mathbf{p}_{k_i}^TB_{k_i}\mathbf{p}_{k_i}
+o\left(\frac{\alpha_{k_i}^2}{\gamma_1^2}\right)&>0.
\end{align*}
The above inequality is equivalent to
\begin{equation}\label{L4.51}
(1-\eta_1)\mathbf{g}_{k_i}^T\mathbf{p}_{k_i}+\frac{\alpha_{k_i}}{2\gamma_1}\left((1-\eta_1)\mathbf{p}_{k_i}^TB_{k_i}\mathbf{p}_{k_i}\right)-\frac{\alpha_{k_i}^2\eta_1\sigma_k}{3\gamma_1^2}\Vert \mathbf{p}_{k_i}\Vert^3+o\left(\frac{\alpha_{k_i}}{\gamma_1}\right)>0.
\end{equation}
Since $\{\mathbf{x}_{k_i}\}$ is on the the unit sphere, there exists a subsequence of $\{\mathbf{x}_{k_i}\}$ that converges to a limit point $\tilde{x}.$ Without confusion, we use $\{\mathbf{x}_{k_i}\}$ to refer to the subsequence whereafter in this proof. Therefore, we have $\mathbf{g}_{k_i}\to \mathbf{g}_{\infty}$,  $B_{k_i} \to B_{\infty}$ and $\mathbf{p}_{k_i}\to \mathbf{p}_{\infty}$. By taking $i\rightarrow\infty,$ we get
$\mathbf{g}_{\infty}^T \mathbf{p}_{\infty}\geq 0$ from \eqref{L4.51} and the hypothesis that $ \alpha_{k_i}\rightarrow0$. On the other hand, in the iteration process
 $
\mathbf{g}_k^T\mathbf{p}_k\leq0.
$
These two inequalities indicate that $$\mathbf{g}_{\infty}^T \mathbf{p}_{\infty}= 0.$$
Thus when $i$ is large enough, from \eqref{L4.51} we get
\begin{equation}\label{L4.52}
\frac{1-\eta_1}{2}\mathbf{p}_{k_i}^TB_{k_i}\mathbf{p}_{k_i}-\frac{\alpha_{k_i}\eta_1\sigma_k}{3\gamma_1}\Vert \mathbf{p}_{k_i}\Vert^3>0.
\end{equation}
Taking $i\rightarrow\infty,$ the above inequality becomes
\begin{equation}\label{L4.53}
 \mathbf{p}_{\infty}^TB_{\infty}\mathbf{p}_{\infty}>0.
\end{equation}

Next, we show
\begin{equation}\label{pBp<0}
 \mathbf{p}_{\infty}^TB_{\infty}\mathbf{p}_{\infty}<0
\end{equation}
based on \eqref{CartisLm1}.  In fact \[
f(\mathbf{x}_k)-m_k(\mathbf{p}_k)\geq 0
\] in Lemma \ref{CartisLm} means that
\begin{equation} \label{L4.54}
-\mathbf{g}_{k_i}^T\mathbf{p}_{k_i}-\frac{1}{2}\mathbf{p}_{k_i}B_{k_i}\mathbf{p}_{k_i}-\frac{1}{3}\sigma_{k_i}\Vert\mathbf{p}_{k_i}\Vert^3 \geq 0.
\end{equation}
By taking $i\to \infty$, we obtain $\mathbf{p}_{\infty}^TB_{\infty}\mathbf{p}_{\infty}\leq -\frac{2}{3}\sigma_{\infty}\Vert\mathbf{p}_{\infty}\Vert^3\leq0,$ which contradicts with the inequality \eqref{L4.53}.

Hence,  $\{\alpha_{k}\}$ is bounded above $0$ and there exists a positive number $\alpha_{\min}$ such that $\alpha_k\geq \alpha_{\min}$ for all $k$.
\end{proof}
Next  We prove $\liminf\limits_{k\rightarrow\infty}\|\mathbf{g}_k\|=0$ by reductio ad absurdum. Assume there exists a constant $ \epsilon >0$ such that $\|\mathbf{g}_k\|\geq \epsilon $ for any $k.$  First we illustrate that $\sqrt{\frac{\Vert \mathbf{g}_k\Vert}{\sigma_k}}$ is a convergent series in Lemma \ref{Lemma3}. Then we show that the parameter $\sigma_k$  is monotonically non-increasing when $k$ is large enough in Lemma  \ref{Lemma4}.  Finally, contradiction emerges on the basis of these two Lemmas.

\begin{Lemma}\label{Lemma3}
Suppose $\mathbf{x}_k,$ $\mathbf{g}_k$ and $\sigma_k$ are produced from Algorithm \ref{algorithm1}. Then we have
\begin{equation}\label{ConvgSeries}
  \sum\limits_{k=1}^{+\infty}\sqrt{\frac{\Vert \mathbf{g}_k\Vert}{\sigma_k}}<+\infty
\end{equation}
and
\begin{equation}\label{p_kUpBound}
 \Vert \mathbf{p}_k \Vert \leq 3\sqrt{\frac{\Vert \mathbf{g}_k\Vert}{\sigma_k}}
\end{equation}
under the assumption $\|\mathbf{g}_k\|\geq \epsilon.$
\end{Lemma}
\begin{proof}
From the definition of $m_k$ in \eqref{sub} we have
\begin{eqnarray}
  m_k(\mathbf{p}_k)-f(\mathbf{x}_k)&=& \mathbf{g}_k^T\mathbf{p}_k+\frac{1}{2}\mathbf{p}_k^TB_k\mathbf{p}_k+\frac{1}{3}\sigma_k\|\mathbf{p}_k\|^3 \label{L4.34-1}\\
  m_k(\alpha\mathbf{p}_k)-f(\mathbf{x}_k) &=&  \alpha\mathbf{g}_k^T\mathbf{p}_k+\frac{1}{2}\alpha^2\mathbf{p}_k^TB_k\mathbf{p}_k+\frac{1}{3}\sigma_k\alpha^3\|\mathbf{p}_k\|^3 \label{L4.34-2}.
\end{eqnarray}
The process of generating $\mathbf{p}_k $ in Algorithm \ref{algorithm1} means that
$$ m_k(\mathbf{p}_k)-f(\mathbf{x}_k)\leq m_k(\mathbf{p}_k^C)-f(\mathbf{x}_k)\leq0.$$
Therefore $\mathbf{g}_k^T\mathbf{p}_k+\frac{1}{2}\mathbf{p}_k^TB_k\mathbf{p}_k\leq 0$ and for any $\alpha_k\in(0,1]$ we have
$$\alpha_k^3(\mathbf{g}_k^T\mathbf{p}_k+\frac{1}{2}\mathbf{p}_k^TB_k\mathbf{p}_k)
\geq\alpha_k^2(\mathbf{g}_k^T\mathbf{p}_k+\frac{1}{2}\mathbf{p}_k^TB_k\mathbf{p}_k)
\geq\alpha_k\mathbf{g}_k^T\mathbf{p}_k+\frac{1}{2}\alpha_k^2\mathbf{p}_k^TB_k\mathbf{p}_k.$$
Combining \eqref{L4.34-1}  with \eqref{L4.34-2}, we get
$$f(\mathbf{x}_k)-m_k(\alpha_k\mathbf{p}_k)\geq\alpha_k^3(f(\mathbf{x}_k)-m_k(\mathbf{p}_k)) .$$
Therefore
\begin{align} \label{L4.30}
& \quad f(\mathbf{x}_k)-f(\mathbf{x}_{k+1}) \geq \eta_1 [f(\mathbf{x}_k)-m_k(\alpha_k\mathbf{p}_k)]\nonumber \geq \alpha_k^3\eta_1 [f(\mathbf{x}_k)-m_k(\mathbf{p}_k)] \nonumber\\
                                    &\geq \frac{\alpha_{\min}^3\eta_1 \|\mathbf{g}_k\| }{6\sqrt{2}}\min \ \left [\frac{\|\mathbf{g}_k\|}{1+\|B_k\|},\frac{1}{2}\sqrt{\frac{\Vert \mathbf{g}_k\Vert}{\sigma_k}}\right] \quad \left(\text{Due to }  \ref{CartisLm1} \ \text{and Lemma} \ \ref{alphakBound}\right) \nonumber\\
                                    & \geq \frac{\alpha_{\min}^3\eta_1 \epsilon }{6\sqrt{2}}\min \ \left [\frac{\epsilon}{1+M},\frac{1}{2}\sqrt{\frac{\Vert \mathbf{g}_k\Vert}{\sigma_k}}\right] \quad \quad \left(\text{Because} \ \|\mathbf{g}_k\|\geq \epsilon \ \text{and} \ \|B_k\|\leq M\right)
\end{align}

Since $m_k(\mathbf{p}_k)-f(\mathbf{x}_k)\leq0,$ the sequence $\{f(\mathbf{x}_k)\}$ is monotonically decreasing. Moreover  $\{f(\mathbf{x}_k)\}$ is bounded below, and therefore it is convergent.
Suppose $\lim\limits_{k\rightarrow+\infty} f(\mathbf{x}_k)=f^*.$
Hence, when $k$ is large enough $f(\mathbf{x}_k)-f(\mathbf{x}_{k+1})$ tends to zero. Based on \eqref{L4.30} we obtain
\begin{equation}\label{L4.31}
  \sqrt{\frac{\Vert \mathbf{g}_k\Vert}{\sigma_k}}\rightarrow 0 \quad \text{as} \ k\rightarrow\infty
\end{equation}
  from the last inequality and
 \begin{align}\label{fk0}
 \sum_{k=k_0}^{+\infty}[f(\mathbf{x}_k)-f(\mathbf{x}_{k+1})]&=f(\mathbf{x}_{k_0})-f^*           \nonumber \\
                                              & \geq \sum_{k=k_0}^{+\infty} \frac{\alpha_{\min}^3\eta_1 \epsilon }{6\sqrt{2}}\min \ \left [\frac{\epsilon}{1+M},\frac{1}{2}\sqrt{\frac{\Vert \mathbf{g}_k\Vert}{\sigma_k}}\right]\nonumber.
\end{align}
Since $f(\mathbf{x}_k)$ is bounded, the conclusion in \eqref{ConvgSeries} holds immediately.

Since $\|\mathbf{g}_k\|\geq \epsilon,$ then we get $\sqrt{\sigma_k \|\mathbf{g}_k\|}\geq \epsilon \sqrt{\frac{\sigma_k}{\|\mathbf{g}_k\|}}\rightarrow+\infty$
from \eqref{L4.31}. With the help of  \eqref{CartisLm2}, the inequality
\begin{equation}\label{L4.32}
  \Vert \mathbf{p}_k \Vert \leq 3\sqrt{\frac{\Vert \mathbf{g}_k\Vert}{\sigma_k}}
\end{equation}
holds for $k$ is large enough.
\end{proof}

Denote an index set $\mathcal{U}:=\{k:\Vert \mathbf{g}_k \Vert \geq \epsilon \ \text{and} \ \sqrt {\frac{\Vert \mathbf{g}_k\Vert}{\sigma_k}}\to 0 \ \text{as k } \to \infty \}.$ We get the following results about  $\sigma_k$ when $k\in\mathcal{U}$ is large enough.
\begin{Lemma}\label{Lemma4}
When $k \in \mathcal{U}$  and $k \to \infty$,
 for  each iteration $k \in \mathcal{U}$ that is sufficiently large, we obtain
\begin{equation}
\sigma_{k+1} \leq \sigma_k. \label{sigma}
\end{equation}

\end{Lemma}
\begin{proof}
Let $\alpha=1.$ With the help of Lemma \ref{Upperbound}, it can be deduced from the Taylor expansion of $f(\mathbf{x}_k^+(1))$ at the point $\mathbf{x}_k$ that
\[
f(\mathbf{x}_k^+(1))\leq f_k+ \mathbf{g}_k^T(\mathbf{x}_k^+(1)-\mathbf{x}_k)+\frac{M}{2}\Vert \mathbf{x}_k^+(1)-\mathbf{x}_k \Vert^2.
\]
Then
\begin{equation}\label{L4.11}
f(\mathbf{x}_k^+(1))-m_k(\mathbf{p}_k)\leq \mathbf{g}_k^T(\mathbf{x}_k^+(1)-\mathbf{x}_k-\mathbf{p}_k)+\frac{M}{2}(\Vert \mathbf{x}_k^+(1)-\mathbf{x}_k \Vert^2+\Vert\mathbf{p}_k\Vert^2)-\frac{\sigma_k}{3}\Vert\mathbf{p}_k\Vert^3.
\end{equation}
From \eqref{x_kg_k=0} and \eqref{xk+}, we obtain
\begin{equation}\label{L4.12}
 | \mathbf{g}_k^T(\mathbf{x}_k^+(1)-\mathbf{x}_k-\mathbf{p}_k)|= |\mathbf{g}_k^T \mathbf{p}_k | \left| \frac{\Vert\mathbf{p}_k\Vert^2-(\mathbf{p}_k^T\mathbf{x}_k)^2}{4+\Vert\mathbf{p}_k\Vert^2-(\mathbf{p}_k^T\mathbf{x}_k)^2}\right| \leq \frac{ |\mathbf{g}_k^T\mathbf{p}_k|}{4}\Vert\mathbf{p}_k\Vert^2 \leq \frac{M}{4}\Vert \mathbf{p}_k\Vert^3.
\end{equation}
Because
$\mathbf{x}_k^+(\alpha_k)$ and $\mathbf{x}_k$  belong to $\mathbb{S}^{n-1}$ and $\alpha \in (0,1]$, by (\ref{xk+})  we have
\begin{equation}\label{L4.33}
\Vert\mathbf{x}_k^+(\alpha)-\mathbf{x}_k\Vert^2
=2-2(\mathbf{x}_k^+(\alpha))^T\mathbf{x}_k
=\frac{4\alpha^2(\Vert\mathbf{p}_k\Vert^2-(\mathbf{p}_k^T\mathbf{x}_k)^2)}{4+\alpha^2\Vert\mathbf{p}_k\Vert^2-\alpha^2(\mathbf{p}_k^T\mathbf{x}_k)^2}
\leq \Vert\mathbf{p}_k\Vert^2.
\end{equation}
Then we get
\begin{equation}\label{L4.13}
\Vert\mathbf{x}_k^+(1)-\mathbf{x}_k\Vert^2 \leq \Vert\mathbf{p}_k\Vert^2.
\end{equation}
Combining \eqref{L4.11}, \eqref{L4.12} and \eqref{L4.13}, we have
\[
f(\mathbf{x}_k^+(1))-m_k(\mathbf{p}_k)\leq\frac{M}{4}\Vert\mathbf{p}_k\Vert^3+M\Vert\mathbf{p}_k\Vert^2-\frac{\sigma_k}{3}\Vert\mathbf{p}_k\Vert^3.
\]
Substituting \eqref{L4.32} into the above inequality, we obtain
\begin{equation}\label{f1m}
f(\mathbf{x}_k^+(1))-m_k(\mathbf{p}_k) \leq \left[9M \sqrt{\frac{\Vert \mathbf{g}_k\Vert}{\sigma_k}}+(\frac{27M}{4}-9\sigma_k) \frac{\Vert \mathbf{g}_k\Vert}{\sigma_k}  \right]\sqrt{\frac{\Vert \mathbf{g}_k\Vert}{\sigma_k}},
\end{equation}
for all $k \in \mathcal{U}$ sufficiently large.

On the other hand, since $\|\mathbf{g}_k\| \geq \epsilon$ and  $\|B_k\|\leq M,$ we have
\[
f(\mathbf{x}_k)-m_k(\mathbf{p}_k)\geq \frac{\epsilon}{6\sqrt{2}}\min \left[\frac{\epsilon}{1+M},\frac{1}{2}\sqrt{\frac{\Vert \mathbf{g}_k\Vert}{\sigma_k}}\right], \qquad \text{for all } k \in \mathcal{U},
\]
from \eqref{CartisLm1}.
Hence, for  all $ k \in \mathcal{U}$  sufficiently large, the above  inequality means
\begin{equation}\label{fm}
f(\mathbf{x}_k)-m_k(\mathbf{p}_k)\geq \frac{\epsilon}{12\sqrt{2}}\sqrt{\frac{\Vert \mathbf{g}_k\Vert}{\sigma_k}}.
\end{equation}
By the rule of updating $\sigma_k$ in Algorithm \ref{algorithm1}, we have
\begin{eqnarray*}
  && k\text{th iteration is very successful} \\
  \iff && \rho_k>\eta_2 \ \text{and }\ \alpha_k=1  \\
  \iff && r_k {:=} f(\mathbf{x}_k^+(1))-f(\mathbf{x}_k)-\eta_2[m_k(\mathbf{p}_k)-f(\mathbf{x}_k)]<0.
\end{eqnarray*}
The formula of $r_k$ can  be equivalently expressed as
\begin{equation}\label{rk}
r_k=f(\mathbf{x}_k^+(1))-m_k(\mathbf{p}_k)+(1-\eta_2)[m_k(\mathbf{p}_k)-f(\mathbf{x}_k)].
\end{equation}
Combining (\ref{f1m}) (\ref{fm}) with (\ref{rk}), we obtain
\begin{equation}\label{rkk}
r_k \leq \sqrt{\frac{\Vert \mathbf{g}_k\Vert}{\sigma_k}}\left[9M \sqrt{\frac{\Vert \mathbf{g}_k\Vert}{\sigma_k}}+(\frac{27M}{4}-9\sigma_k) \frac{\Vert \mathbf{g}_k\Vert}{\sigma_k}-\frac{(1-\eta_2)\epsilon}{12\sqrt{2}}\right],
\end{equation}
for all $k \in \mathcal{U}$ sufficiently large.
Because $\sqrt {\frac{\Vert \mathbf{g}_k\Vert}{\sigma_k}}\to 0 \ \text{as k } \to \infty $, the inequality \eqref{rkk}
indicates $r_k<0$ for all $k \in \mathcal{U}$ sufficiently large. Therefore,  the $k$th iteration is very successful when $k \in \mathcal{U}$ is sufficiently large. Following from the updating procedure of $\sigma_k$  in Algorithm \ref{algorithm1}, the inequality (\ref{sigma}) holds.
\end{proof}

According to the termination criteria in Algorithm \ref{algorithm1}, if the algorithm terminates finitely, then $\mathbf{g}(\mathbf{x}_*)=0$  at the end point $\mathbf{x}_*.$ Next we demonstrate that there exists a subsequence of $\{\mathbf{g}_k\}$ converging to zero  when iteration points are infinite.
\begin{Theorem}\label{infgktend0}
Suppose the infinite sequence $\{\mathbf{x}_k\}$ is produced by Algorithm \ref{algorithm1}. Then  its gradient sequence   $\{\mathbf{g}_k\}$ satisfies
\begin{equation}\label{infgktend0Eq}
\liminf_{k\to\infty}\Vert \mathbf{g}_k \Vert=0.
\end{equation}
\end{Theorem}
\begin{proof}
Assume $\|\mathbf{g}_k\|\geq \epsilon.$ From \eqref{ConvgSeries} in Lemma \ref{Lemma3}, we have
$$\sigma_k\rightarrow +\infty,$$
which is incompatible with the conclusion
$$\sigma_{k+1}\leq \sigma_k ,\quad  \text{for sufficiently large} \ k$$
in Lemma \ref{Lemma4}. Therefore, the assumption is invalid and \eqref{infgktend0Eq} holds.
\end{proof}

Based on the conclusion in Theorem \ref{infgktend0}, we prove that the whole sequence $\{\mathbf{g}_k\}$ converges to zero.
\begin{Theorem}\label{g_ktend0}
 Suppose the infinite sequence $\{\mathbf{x}_k\}$ is produced by Algorithm \ref{algorithm1}. Then  its gradient sequence   $\{\mathbf{g}_k\}$ satisfies
\begin{equation}\label{g_ktend0E}
\lim_{k\to \infty}\Vert \mathbf{g}_k\Vert=0.
\end{equation}
\end{Theorem}
\begin{proof}
Assume \eqref{g_ktend0E} does not hold and there exists an infinite subsequence of iterations $\{t_i\}$ such that
\begin{equation}\label{T4.20}
\Vert \mathbf{g}_{t_i}\Vert \geq 2 \epsilon, \qquad \text{for some }\epsilon >0\text{ and for all }i.
\end{equation}
On the other hand, Theorem \ref{infgktend0} indicates that there is an infinite subsequence of $\{\mathbf{g}_k\}$ that converges to zero. Therefore, we choose the iteration sequence $\{l_i\}$ such that $l_i$ is the first iteration satisfying
\begin{equation}\label{T4.24}
 \|\mathbf{g}_{l_i}\|\leq \epsilon
\end{equation}
after $t_i,$  which means that for all $i$, we have
\begin{equation}\label{T4.21}
\Vert \mathbf{g}_m \Vert \geq \epsilon, \qquad \text{for all }m\text{ with }t_i\leq m <l_i.
\end{equation}
From \eqref{L4.30}, \eqref{L4.31} and \eqref{L4.32}, we get
\begin{equation}\label{T4.23}
f(\mathbf{x}_m)-f(\mathbf{x}_{m+1})\geq  \frac{\alpha_{\min}^3\eta_1 \epsilon}{36\sqrt{2}}\Vert \mathbf{p}_m \Vert,\quad \text{for all }t_i\leq m<l_i,\ i\text{ sufficiently large.}
\end{equation}
Summing up \eqref{T4.23} from $t_i$ to $l_{i-1},$ we obtain
 \begin{align*}
 \frac{36\sqrt{2}}{\eta_1\epsilon}[f(\mathbf{x}_{t_i})-f(\mathbf{x}_{l_i})]&\geq \sum_{k=t_i}^{l_i-1}\alpha_{\min}^3\Vert \mathbf{p}_k\Vert\\
                & \geq \sum_{k=t_i}^{l_i-1}\alpha_{\min}^3\Vert \mathbf{x}_{k+1}-\mathbf{x}_k\Vert  \quad \left(\text{by} \ \eqref{L4.33}\right)\\
                                                &\geq \alpha_{\min}^3\Vert \mathbf{x}_{t_i}-\mathbf{x}_{l_i}\Vert
 \end{align*}
for all $ i$ sufficiently large.

In the proof process of Lemma \ref{Lemma3}, we have explained that the sequence $\{f(\mathbf{x}_k)\}$ is  convergent. Thus
$\left[f(\mathbf{x}_{t_i})-f(\mathbf{x}_{l_i})\right]\rightarrow 0$ as $i\rightarrow +\infty.$ Then the sequence $\|\mathbf{x}_{t_i}-\mathbf{x}_{l_i}\|$ converges to zero from the above inequalities. Since $\nabla^2 f( \mathbf{x})$ is bounded, $ \mathbf{g(x)}$ is uniformly continuous, which implies that $\|\mathbf{g}_{t_i}-\mathbf{g}_{l_i}\|$ tends to zero as $\|\mathbf{x}_{t_i}-\mathbf{x}_{l_i}\|$ converges to zero. However, from \eqref{T4.20} and \eqref{T4.24} we have
$$\Vert \mathbf{g}_{t_i}-\mathbf{g}_{l_i}\Vert \geq \Vert\mathbf{g}_{t_i}\Vert-\Vert\mathbf{g}_{l_i}\Vert \geq \epsilon,$$
which produces a contradiction. Hence, the assumption \eqref{T4.20} does not hold. The proof is completed.
\end{proof}
\section{Fast computation skill for $\mathcal{T}\mathbf{x}^{r-2}$ }

The eigenvalues of a symmetric tensor play an important role in spectral hypergraph theory. In the  process of tensor eigenvalue computation, the operator $\mathcal{T}\mathbf{x}^{r-2}$ is frequently invoked. In this section, we introduce a fast computation skill FCS for tensor-vector product $\mathcal{T}\mathbf{x}^{r-2}$, in which the tensor $\mathcal{T}$ is arisen from a hypergraph.

\subsection{The basics of hypergraphs}
\begin{Definition}[Hypergraph] A hypergraph is defined as $G=(V,E)$, where $V=\{1,2,\cdots,n\}$ is the vertex set and $E=\{e_1,e_2,\cdots,e_m\}$ is the edge set for $e_p\subset V,p=1,2,\cdots,m.$ If $|e_p|=r\geq 2$ for $p=1,2,\cdots,m$ and $e_i \neq e_j$ when $i\neq j$, we call $G$ an $r$-uniform hypergraph or an $r$-graph. If $r=2$, $G$ is an ordinary graph.
\end{Definition}

For each vertex $i\in V$, the degree of $i$ is defined as
\[
d(i)=|\{e_p:i\in e_p,e_p\in E\}|.
\]
\begin{figure}[ht]
\centering
\includegraphics[height=3cm]{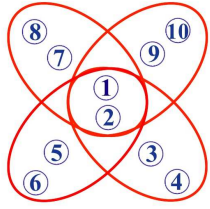}
\caption{A 4-uniform flower hypergraph.}
\label{figf}
\end{figure}

For instance, the flower hypergraph illustrated in Figure \ref{figf} is a 4-uniform hypergraph. There are ten vertices $V=\{1,2,\cdots,10\}$ and four edges $E=\{e_1=\{1,2,3,4\},e_2=\{1,2,5,6\},e_3=\{1,2,7,8\},e_4=\{1,2,9,10\}\}$ in this hypergraph. Its all edges share two common vertices, and the degree of each vertex is $d(1)=d(2)=4$ and $d(i)=1$ for $i=3,4,\cdots,10.$

\begin{Definition}[adjacency tensor\cite{cooper2012spectra}] For an $r$-graph $G=(V,E)$ with $n$ vertices, the adjacency tensor $\mathcal{A}=(a_{i_1\cdots i_r})$ of $G$ is an $r$th order $n$-dimensional symmetric tensor with  entries
\[
a_{i_1\cdots i_r}=\begin{cases}
\frac{1}{(r-1)!}  & \text{if } \{i_1,\cdots,i_r\} \in E,\\
0                 & \text{otherwise.}
\end{cases}
\]
\end{Definition}

\begin{Definition}[Laplacian tensor and signless Laplacian tensor\cite{2013H}] For an $r$-graph $G=(V,E)$ with $n$ vertices, the degree tensor $\mathcal{D}$ is defined as an $r$th order $n$-dimensional diagonal tensor whose $i$th diagonal element is $d(i)$. Then the Laplacian tensor $\mathcal{L}$ of $G$ is defined as
\[
\mathcal{L}=\mathcal{D}-\mathcal{A},
\]
and the signless Laplacian tensor $\mathcal{Q}$ of $G$ is defined as
\[\mathcal{Q}=\mathcal{D}+\mathcal{A}.
\]

\label{ALQ}
\end{Definition}
\subsection{Calculation of $\mathcal{T}\mathbf{x}^{r-2}$}
In order to compute the extreme Z-eigenvalue and H-eigenvalue of the adjacency tensor $\mathcal{A}$, the Laplacian tensor $\mathcal{L}$ and the signless Laplacian tensor $\mathcal{Q}$,  we need to compute products $\mathcal{T}\mathbf{x}^r$, $\mathcal{T}\mathbf{x}^{r-1}$  and $\mathcal{T}\mathbf{x}^{r-2}$ when $\mathcal{T}=\mathcal{A},\ \mathcal{L}\text{ and }\mathcal{Q}$. Chang et al.\cite{chang2016computing} provided an economical way to store a uniform hypergraph. Based on the economical storage, a fast method for computing products $\mathcal{T}\mathbf{x}^r$ and $\mathcal{T}\mathbf{x}^{r-1}$ was established. In this section, we introduce a fast tensor-vector product skill for $\mathcal{T}\mathbf{x}^{r-2}$, for $\mathcal{T}=\mathcal{A},\ \mathcal{L}\text{ and }\mathcal{Q}$.

Let $G=(V,E)$ be an $r$-uniform hypergraph with $n$ vertices and $m$ edges. We store $G$ as an $m$-by-$r$ matrix $G_m$ whose each row represents an edge of $G$, and the entries in each row are vertices.

For example, we consider the 4-uniform flower hypergraph $G$ which is shown in Figure \ref{figf}. Then $G$ can be stored by a $4-$by$-4$ matrix
\begin{equation} G_m= \begin{bmatrix}  \label{GM}
1 & 2 & 3 & 4 \\
1 & 2 & 5 & 6 \\
1 & 2 & 7 & 8 \\
1 & 2 & 9 & 10
\end{bmatrix}\in \mathbb{R}^{m\times r}.\end{equation}

 Consider the degree tensor $\mathcal{D}$, whose $i$th diagonal element is the degree $d(i)$,
\[
d(i)=\sum_{j=1}^r\sum_{l=1}^m\delta (i,(G_m)_{lj}), \ \ \text{for }i=1,2,\cdots,n.
\]
Here $\delta(\cdot,\cdot)$ is the kronecker notation, i.e.
\[
\delta(i,j)=\begin{cases}
1, & \text{if }i=j,\\
0, & \text{if }i\ne j.
\end{cases}
\]
Then for any vector $\mathbf{x}\in \mathbb{R}^n$,
\[
\mathcal{D}\mathbf{x}^{r-2}=diag\{d(1)x_1^{r-2},\cdots,d(n)x_n^{r-2}\}.
\]

In order to compute $\mathcal{A}\mathbf{x}^{r-2}$, we construct another matrix $X_m=[x_{(G_m)_{lj}}]$ whose size is the same as that of $G_m$. If the $(l,j)$-th element of $G_m$ is the vertex $i$, then the $(l,j)$-th element of $X_m$ is $x_i$. For example, the matrix $X_m$ corresponding to the matrix $G_m$ in (\ref{GM}) is
\[
 X_m= \begin{bmatrix}
x_1 & x_2 & x_3 & x_4 \\
x_1 & x_2 & x_5 & x_6 \\
x_1 & x_2 & x_7 & x_8 \\
x_1 & x_2 & x_9 & x_{10}
\end{bmatrix}.
\]
Based on the matrix $X_m$, we rewrite the calculation formula of $(\mathcal{A}\mathbf{x}^{r-2})_{ij}$ as
\[(\mathcal{A}\mathbf{x}^{r-2})_{ij}=\sum_{s=1}^m\sum_{\substack {l,k=1\\ l\ne k}}^r\left(\delta(i,(G_m)_{sl})\delta(j,(G_m)_{sk})\prod_{\substack{t\ne l\\t\ne k}}(X_m)_{st}\right),\text{ for }i,j=1,2,\cdots,n.
\]

\begin{figure*}

\lstset{language=Matlab}
\begin{lstlisting}

% Store a 4-uniform flower hypergraph
Gm = [1,2,3,4; 1,2,5,6; 1,2,7,8; 1,2,9,10];

% Calculate the degree vector
[m,r]  = size(Gm);
n      = max(Gm(:));
Md     = sparse(Gm(:),(1:m*r)',ones(m*r,1),n,m*r);
degree = full(sum(Md,2));

% Compute Dx^(r-2)
Dx2  = diag(degree.*(x.^(r-2)));

% Compute Ax^(r-2)
Xm   = reshape(x(Gm(:)),[m,k]);
Axk2 = zeros(n,n);
for i = 1:r-1
   for j = i+1:r
       A = Xm;
       A(:,[i,j]) = [];
       B = prod(A,2);
       Axk2 = Axk2+sparse(Gm(:,i),Gm(:,j),B,n,n);
   end
end
Axk2 = (Axk2+Axk2').*factorial(r-2)/factorial(r-1);

% Compute Lx^(r-2) and Qx^(r-2)
Lxk2 = Dxk2-Axk2;
Qxk2 = Dxk2+Axk2;

\end{lstlisting}
\caption{Code for computing $\mathcal{T}\mathbf{x}^{r-2}$}
\label{code}
\end{figure*}
We show the code for computing $\mathcal{T}\mathbf{x}^{r-2}$ in the figure \ref{code}.

\subsection{Test of the fast computation skill}
In this section, we compare the fast computation skill FCS with the traditional algorithm for computing $\mathcal{A}\mathbf{x}^{r-2}$. Here $\mathcal{A}$ refers to the adjacency tensor  of the 4-uniform flower hypergraph in Figure \ref{figf}. In Tensor Toolbox, the traditional algorithm for computing $\mathcal{A}\mathbf{x}^{r-2}$ is implemented as $\mathbf{ttsv}$. The vector $\mathbf{x} \in \mathbb{R}^n$ is randomly generated. All numerical experiments in this paper are carried by a laptop with i5-10210U CPU at 1.60GHz and 16.0GB of RAM .

\begin{figure*}[htbp] 
\centering
\begin{minipage}[t]{0.54\linewidth} 
\centering
\includegraphics[height=3.8cm]{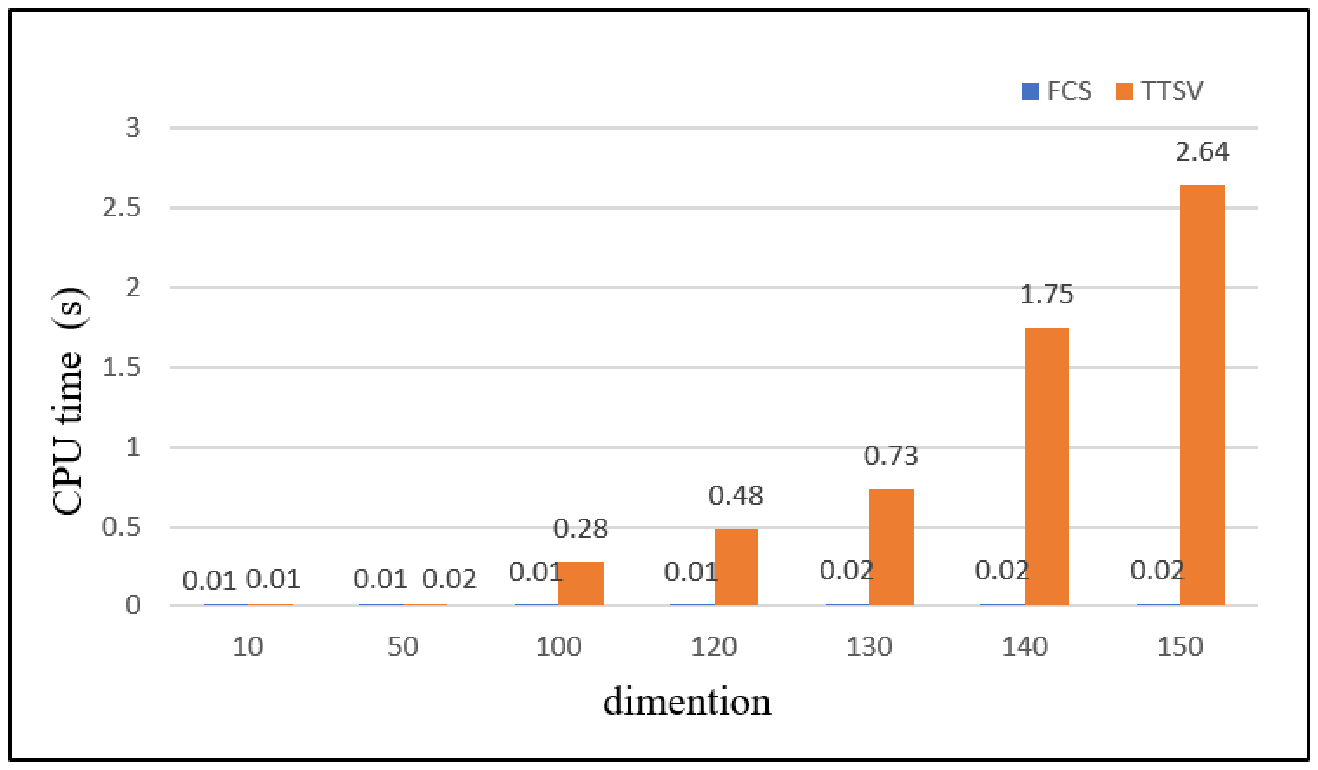} 
\centerline{(a)}
\end{minipage}%
\begin{minipage}[t]{0.52\linewidth}
\centering
\includegraphics[height=3.8cm]{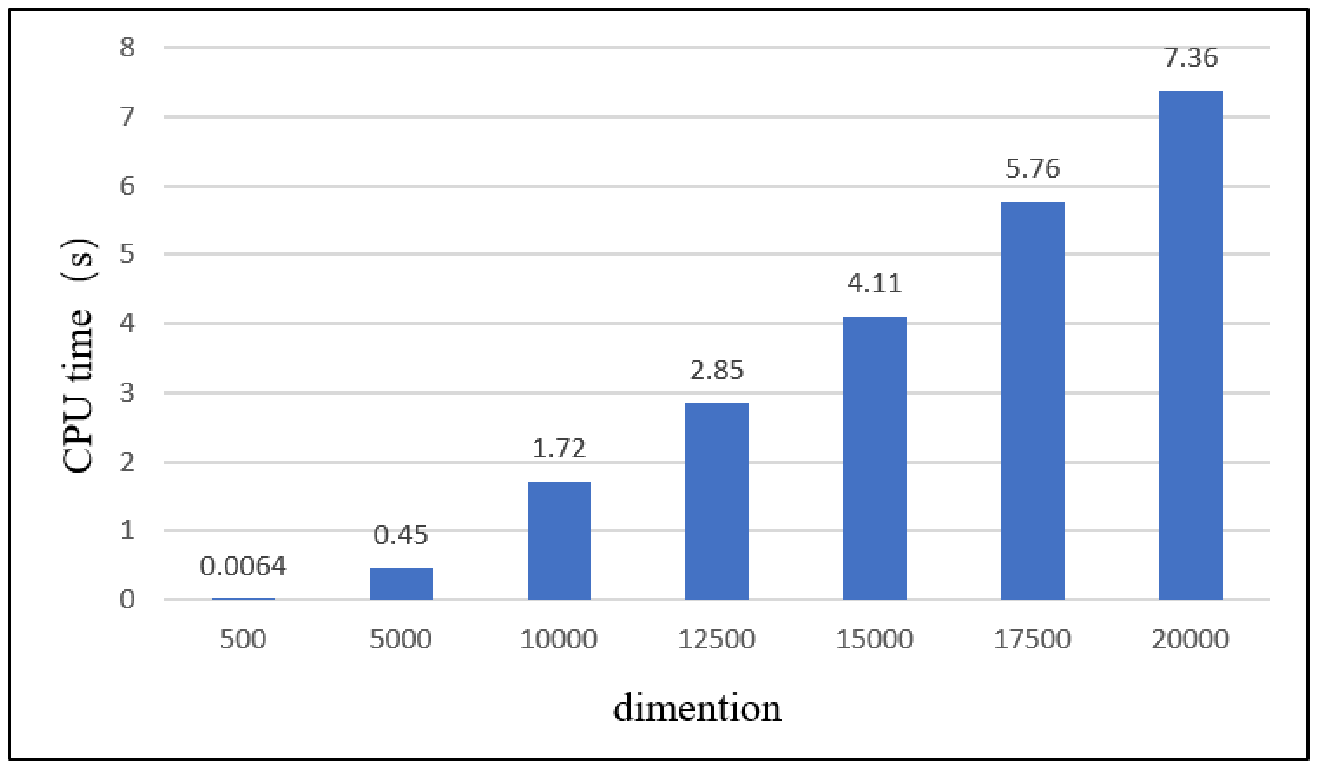}
\centerline{(b)}
\end{minipage}%
\caption{CPU time for computing $\mathcal{A}\mathbf{x}^{r-2}$ with different dimentions $n$.}
\label{FCS}
\end{figure*}

The performance of  computing $\mathcal{A}\mathbf{x}^{r-2}$  by FCS and $\mathbf{ttsv}$ are shown in Figure \ref{FCS}.
 As shown in  Figure \ref{FCS} (a), when $n\leq 150$, it takes about 0.01s for FCS to calculate $\mathcal{A}\mathbf{x}^{r-2}$. With the increase of the dimension $n$, $\mathbf{ ttsv}$  spends much more time  to calculate $\mathcal{A}\mathbf{x}^{r-2}$ than time taken by FCS. Furthermore,  our laptop is incapable of  computing eigenvalues of tensors  with dimension  greater than 150 by $\mathbf{ttsv}$, while FCS can compute tensors of higher dimensions. Figure \ref{FCS} (b) demonstrates the results of CPU time for calculating $\mathcal{A}\mathbf{x}^{r-2}$ with large $n$ by FCS. Up to 20000 dimensions, the calculation by FCS costs less than 8 seconds.

\section{Numerical experiments}

In order to show the efficiency of the proposed adaptive cubic regularization algorithm (ACRCET), we perform some numerical experiments. In this section, three other algorithms are compared with ACRCET .

$\bullet$ An adaptive power method (PM) \cite{kolda2011shifted,kolda2014adaptive}. In Tensor Toolbox, it is implemented as $\mathbf{eig\_sshopm}$ and $\mathbf{eig\_ geap}$ for Z-eigenvalues and H-eigenvalues of even order symmetric tensors, respectively.

$\bullet$ A curvilinear search algorithm (ACSA) \cite{chenq2016computing} which employs the Barzilai-Borwein gradient optimization algorithm for computing tensor eigenvalues. The code was provided by the authors of \cite{chenq2016computing}.

$\bullet$ An unconstrained optimization approach (HUOA) \cite{han2013unconstrained}. The author employs the local optimization solver $\mathbf{fminunc}$ from the Optimization Toolbox to solve an unconstrained optimization model.

The ACRCET algorithm is  implemented with following parameters
\[
\eta_1=0.1, \eta_2=0.5, \gamma_1=0.25, \gamma_2=1.2, \text{and } \gamma_3=2.
\]
We solve the subproblem based on the method given in \cite[Section 6]{cartis2011adaptive}. First we employ the Lanczos method to simplify the subproblem \eqref{sub} as
\begin{equation*}
  m(\mathbf{u}) = f_k + \gamma \mathbf{u}^T \mathbf{e}_1 +\frac{1}{2}\mathbf{u}^T T_k\mathbf{u}+\frac{1}{3}\sigma_k\|\mathbf{u}\|^3,
\end{equation*}
where $\mathbf{e}_1$ is a unit vector and $T$ is a symmetric tridiagonal matrix. Similar to the trust region method, it is proved in
\cite[Theorem3.1]{cartis2011adaptive}  that $\mathbf{u}$ is a global minimizer of the above subproblem if and only if a pair of $(\mathbf{u},\lambda)$ satisfies
\begin{equation} \label{eq6.1}
  (T_k+\lambda I)\mathbf{u}=-\gamma\mathbf{e}_1 \text{ and } \lambda^2=\sigma_k^2\mathbf{u}^T\mathbf{u}
\end{equation}
where $T_k+\lambda I$ is positive semidefinite. The equation system \eqref{eq6.1} is finally solved by Newton's method \cite[Algorithm 6.1]{cartis2011adaptive}.

We compute the extreme H- or Z-eigenvalue of a tensor   by running four algorithms PM, ACSA, HUOA, and ACRCET from 100 random initial points sampled on the unit sphere $\mathbb{S}^{n-1}$. Then we obtain 100 estimated eigenvalues and take the best one as the estimated extreme eigenvalue. In the following experiments, we report the estimated extreme eigenvalue, the total number of iterations (Iter'n) and the total CPU time (Time in seconds) of the 100 runs.

$\mathbf{Example\ 1}$ In \cite{qi2005eigenvalues}, Qi generated a symmetric tensor $\mathcal{A}(\alpha)\in\mathbb{R}^{[4,2]}$  with
\[
a_{1111}=3,a_{2222}=1,\
 a_{1122}=a_{1221}=a_{1212}=a_{2121}=a_{2211}=a_{2112}=\alpha
\]
and other entries being zero. For different values of $\alpha,$ all Z-eigenvalues of $\mathcal{A}(\alpha)$ are analyzed  and provided in \cite{qi2005eigenvalues}. Thus we know the smallest Z-eigenvalue of $\mathcal{A}(\alpha)$ for $\alpha=0,10$ and $100$ are
$$\lambda_{\min}^Z(\mathcal{A}(0))=\frac{3}{4},\lambda_{\min}^Z(\mathcal{A}(10))=\lambda_{\min}^Z(\mathcal{A}(100))=1.$$

We compute the smallest Z-eigenvalues of $\mathcal{A}(\alpha)$ when $\alpha=0,10$ and $100$ by PM, ACSA, HUOA and ACRCET. The numerical results are shown in Table \ref{tab1}. It can be seen that both PM and ACRCET find all the true smallest eigenvalues. When $\alpha=10$ or $100$, ACSA gives one of the Z-eigenvalues of $\mathcal{A}(\alpha)$, but misses the smallest one. When $\alpha=0$ or $10,$ HUOA obtains the smallest Z-eigenvalue  of $\mathcal{A}(\alpha)$ inaccurately. Since second order information from the Hessian and cubic overestimator of the objective function are employed, the proposed ACRCET  performs much better than other algorithms.

\renewcommand{\arraystretch}{2} 
\begin{table*}[]

\scriptsize

\centering

  \fontsize{10}{6}\selectfont

\caption{Results for finding the smallest  Z-eigenvalues of Example 1.}

\label{tab1}
\scalebox{0.78}{
\begin{tabular}{cccc|ccc|ccc}
\hline

\multirow{2}{*}{Algorithms} & \multicolumn{3}{c}{$\alpha=0$} & \multicolumn{3}{c}{$\alpha=10$} & \multicolumn{3}{c}{$\alpha=100$} \\

\cmidrule(r){2-4} \cmidrule(r){5-7}  \cmidrule(r){8-10}

&  $\lambda_{\min}^Z$      &  Iter'n   &   Time(s)
&  $\lambda_{\min}^Z$      &  Iter'n   &   Time(s)
&  $\lambda_{\min}^Z$      &  Iter'n   &   Time(s)\\
\hline
PM        &0.7500     & 813 & 0.26   & 1.0000 & 588 & 0.21 & 1.0000 & 561 & 0.29 \\

ACSA  &0.7500  & 800 & 0.11 & 3.0000(*) & 500 & 0.12 & 3.0000(*) & 500& 0.10 \\

HUOA  &0.7533(*) & 3219 & 1.34& 1.1149(*)& 2893 &1.33 & 1.0000 & 3219 &1.47 \\

ACRCET   &{\bf 0.7500}  &{\bf 200} &{\bf 0.12} & {\bf 1.0000} &{\bf 200} & {\bf 0.13} &{\bf 1.0000}& {\bf 400}&{\bf 0.15} \\

\hline

\end{tabular}
}
\end{table*}

\begin{figure}[ht]
\centering
\includegraphics[height=3cm]{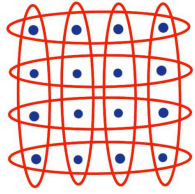}
\caption{A 4-uniform 2-regular hypergraph.}
\label{fig1}
\end{figure}

$\mathbf{Example\ 2\ (\textbf{A 2-Regular}\ Hypergraph)}$. A regular hypergraph is a hypergraph whose vertices have the same degree. For the hypergraph presented in Figure \ref{fig1}, the degree of all its vertices is 2. We call it  a $2$-regular hypergraph and denote it as $G^2_R$.
Define the H-spectral radius of  a tensor $\mathcal{T}$  as the largest modulus of the H-eigenvalues of  $\mathcal{T},$ and we use $\rho(\mathcal{T})$ to stand for the H-spectral radius of  $\mathcal{T}.$
\begin{Proposition}[\cite{qi2017tensor}]\label{d-regular}
 Suppose $G$ is a $d$-regular hypergraph. The H-spectral radius of its signless Laplacian  tensor $\rho(\mathcal{Q}(G))=2d$ and the H-spectral radius of its adjacency tensor $\rho(\mathcal{A}(G))=d.$
\end{Proposition}

We compute the largest H-eigenvalue of the signless Laplacian tensor $\mathcal{Q}(G_R^{2})$ and the smallest H-eigenvalue of the adjacency tensor $\mathcal{A}(G_R^{2}).$ The results are shown in Table \ref{tab2}.  From Proposition \ref{d-regular}, we know $\rho(\mathcal{Q}(G^2_R))=4$ and  $\rho(\mathcal{A}(G^2_R))=2.$ The extreme eigenvalues calculated by PM,  ACSA, HUOA, and ACRCET agree with this conclusion.    It can be seen that ACRCET runs much faster than the other three methods.

\renewcommand{\arraystretch}{1.5} 
\begin{table}[tp]

  \centering
  \fontsize{9}{8}\selectfont
  \begin{threeparttable}
  \caption{Results for finding the largest H-eigenvalue of $\mathcal{Q}(G_R^2)$ and the smallest H-eigenvalue of $\mathcal{A}(G_R^2)$.}
  \label{tab2}
    \begin{tabular}{c|ccc|ccc}
\hline
  Algorithms&$\lambda_{\max}^H(\mathcal{Q}(G_R^2))$&Iter'n&Time(s)&$\lambda_{\min}^H(\mathcal{A}(G_R^2))$&Iter'n&Time(s)\cr
   \hline
    PM&4.0000&4310&2.52&-2.0000&4596&3.06\cr
    ACSA&4.0000&2952&1.97&-2.0000&2876&2.05\cr
    HUOA&4.0000&3304&1.82&-2.0000&3297&1.63\cr
    ACRCET&{\bf 4.0000}&{\bf 616}&{\bf 0.71}&{\bf-2.0000}&{\bf 618}&{\bf 0.60}\cr
    \hline
    \end{tabular}
    \end{threeparttable}
\end{table}

\begin{figure*}[htbp]
\begin{minipage}[t]{0.25\linewidth}
\centering
\includegraphics[height=3cm]{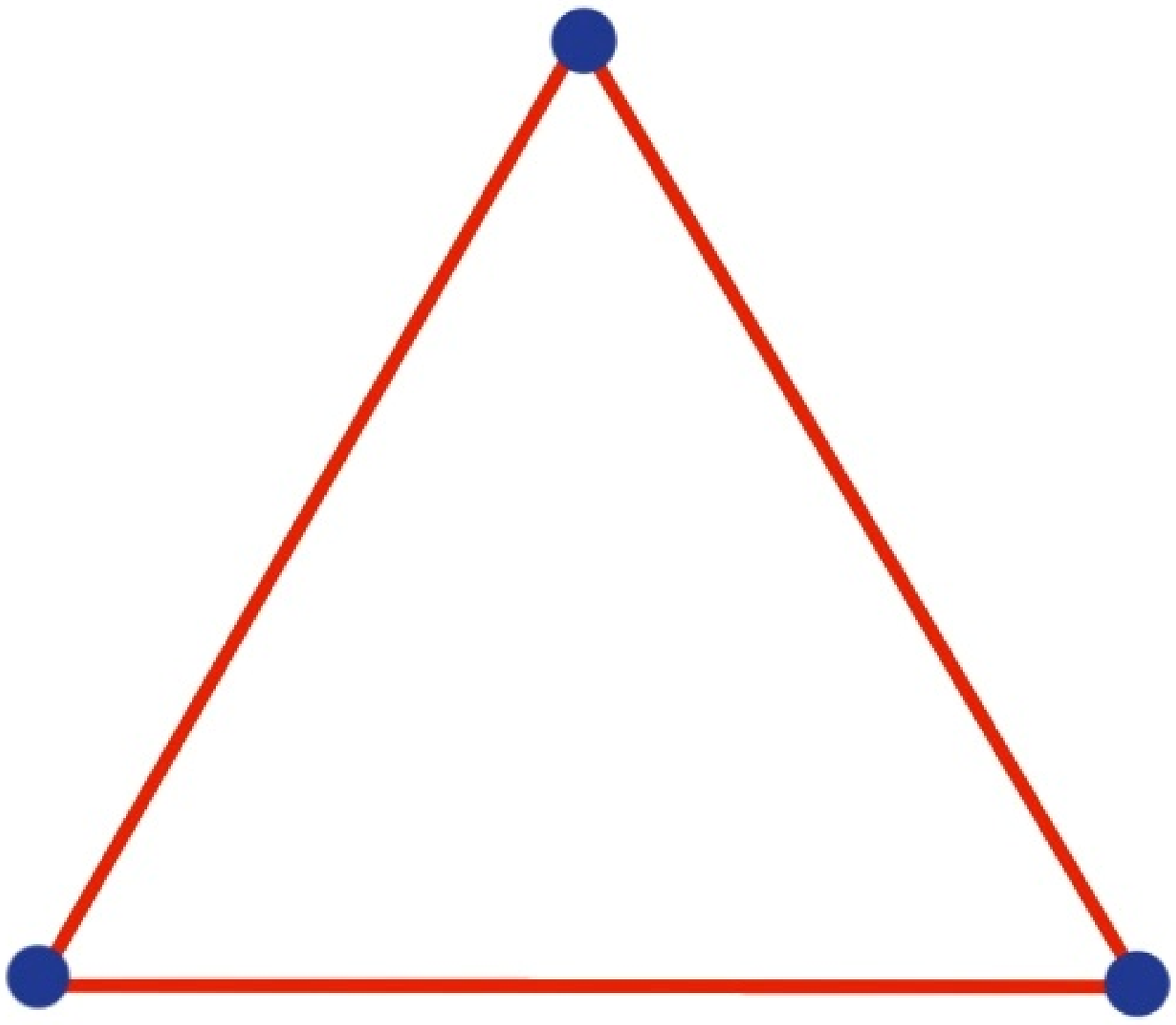}
\centerline{$(a) \ G^{2,3}_L$}
\end{minipage}%
\begin{minipage}[t]{0.25\linewidth}
\centering
\includegraphics[height=3cm]{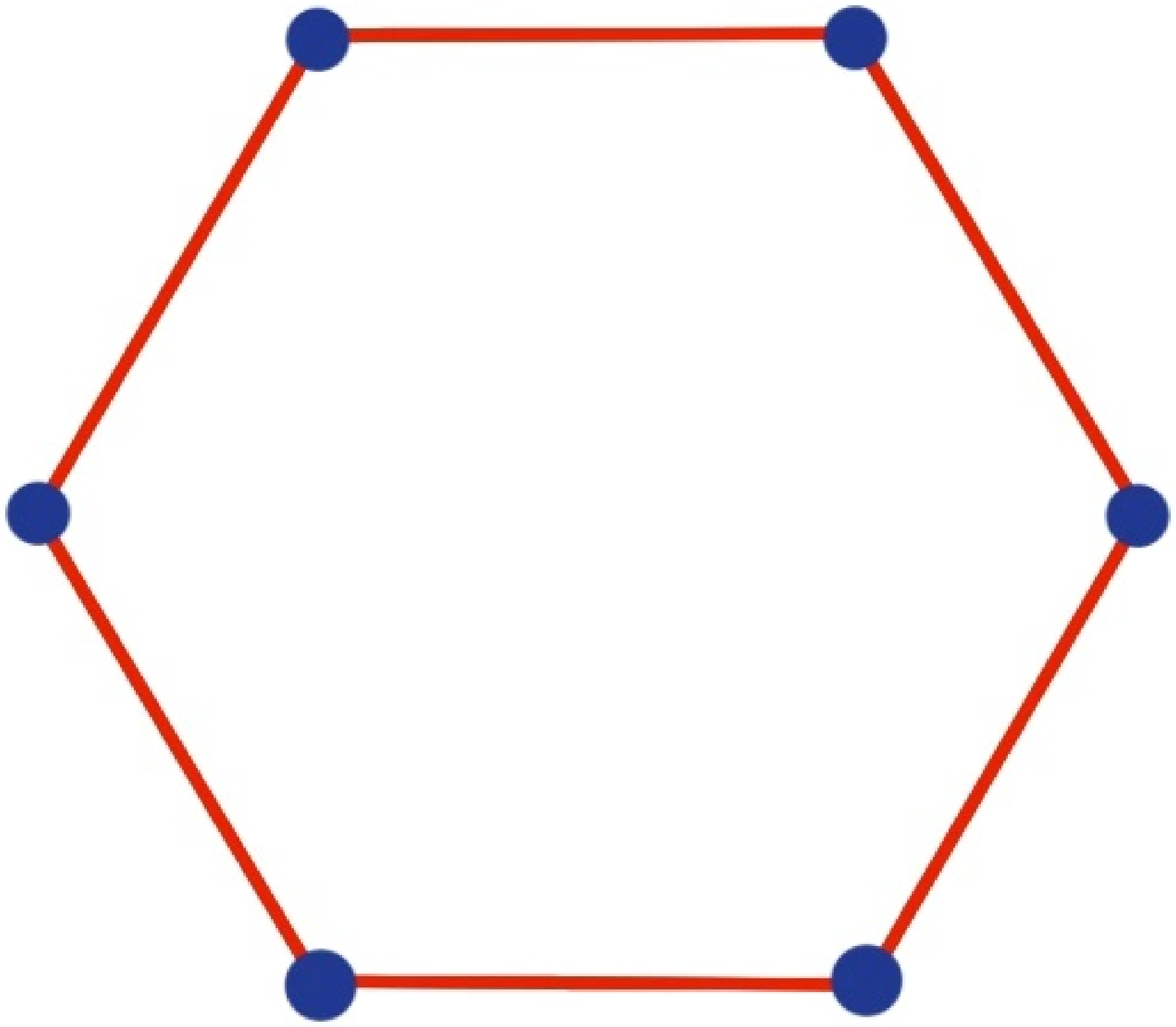}
\centerline{$(b) \ G^{2,6}_L$}
\end{minipage}

\quad

\begin{minipage}[t]{0.25\linewidth}
\centering
\includegraphics[height=3cm]{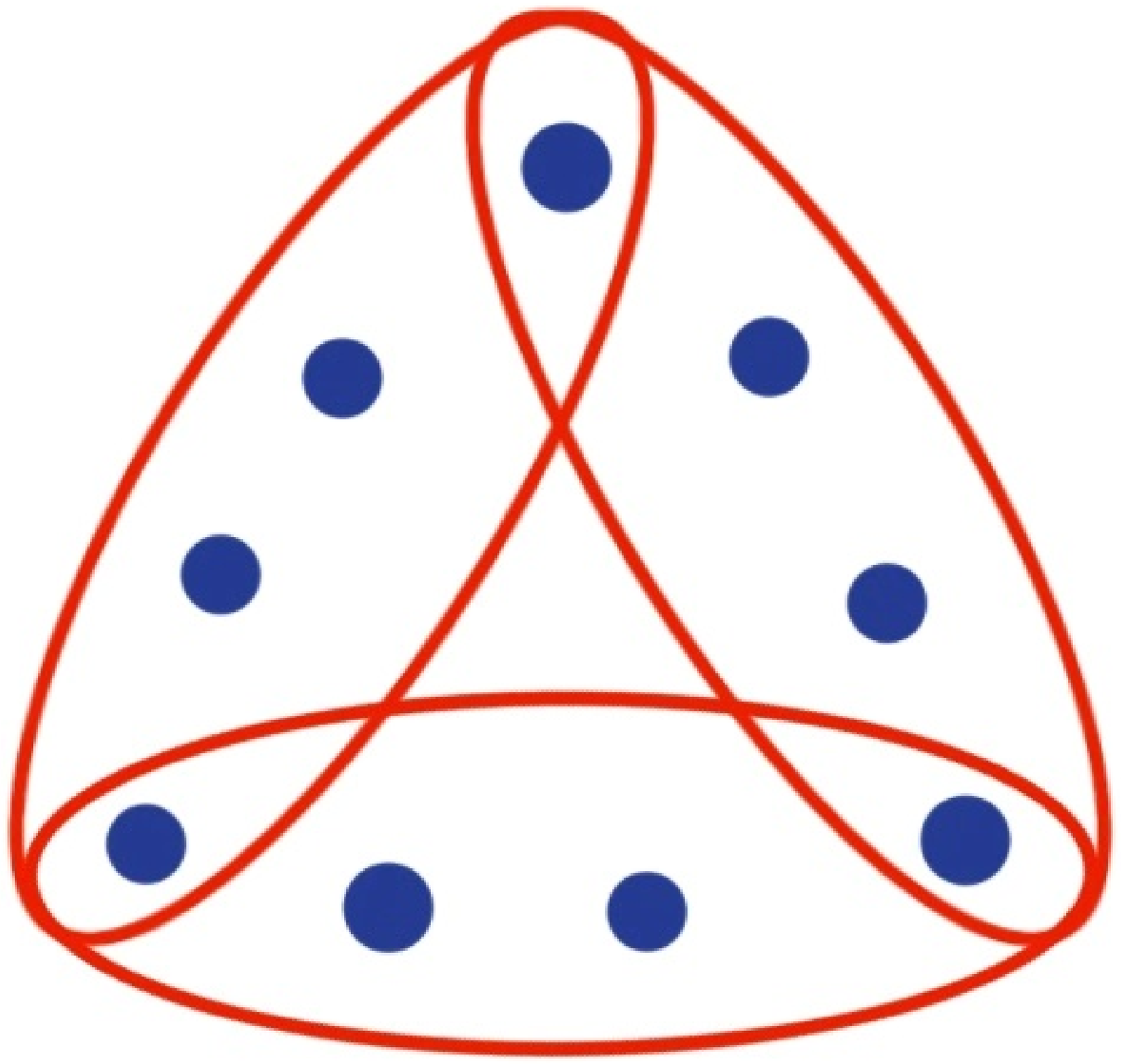}
\centerline{$(c) \ G^{4,3}_L $}
\end{minipage}%
\begin{minipage}[t]{0.25\linewidth}
\centering
\includegraphics[height=3cm]{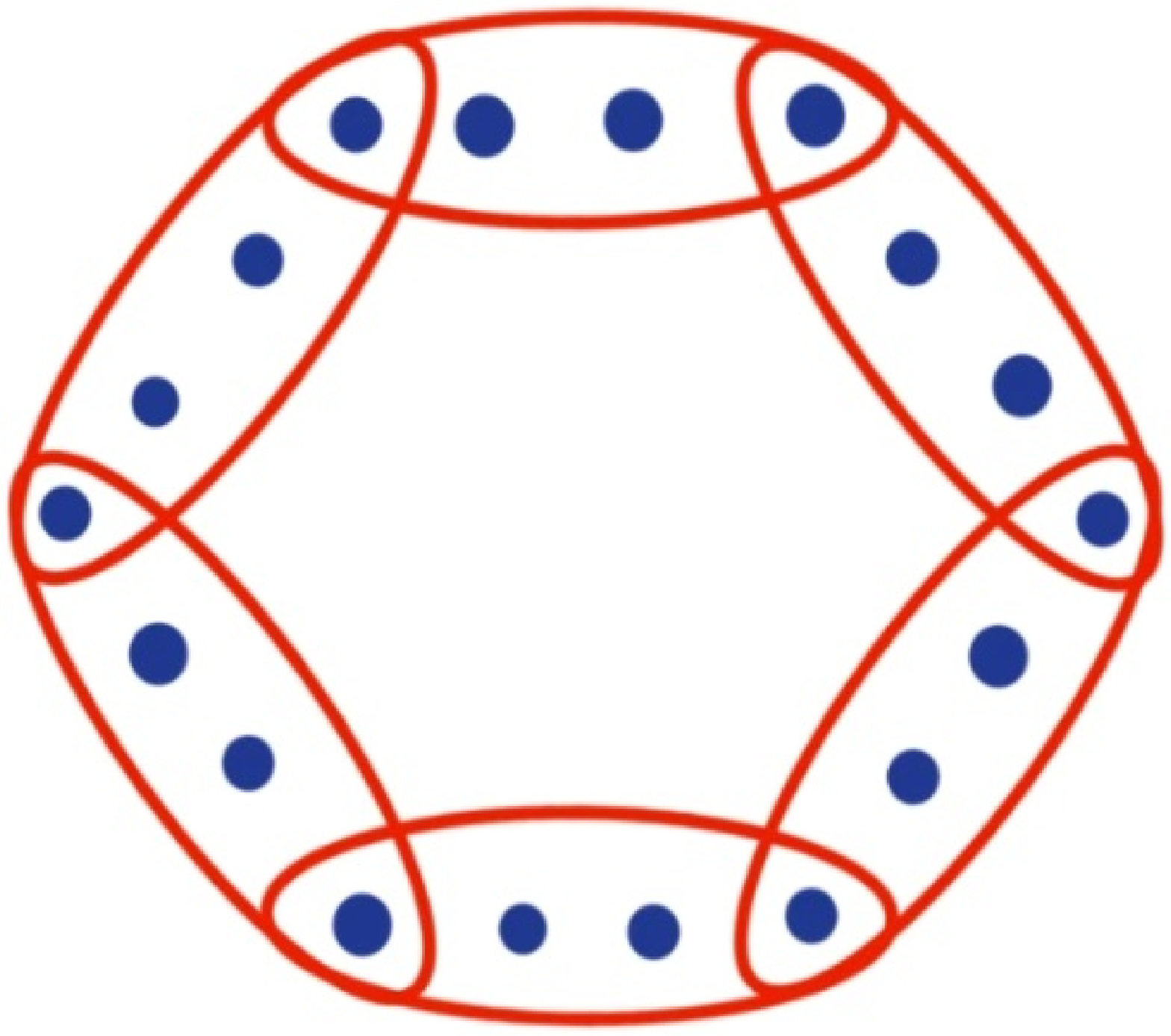}
\centerline{$(d) \ G^{4,6}_L $}
\end{minipage}
\begin{minipage}[t]{0.25\linewidth}
\centering
\includegraphics[height=3cm]{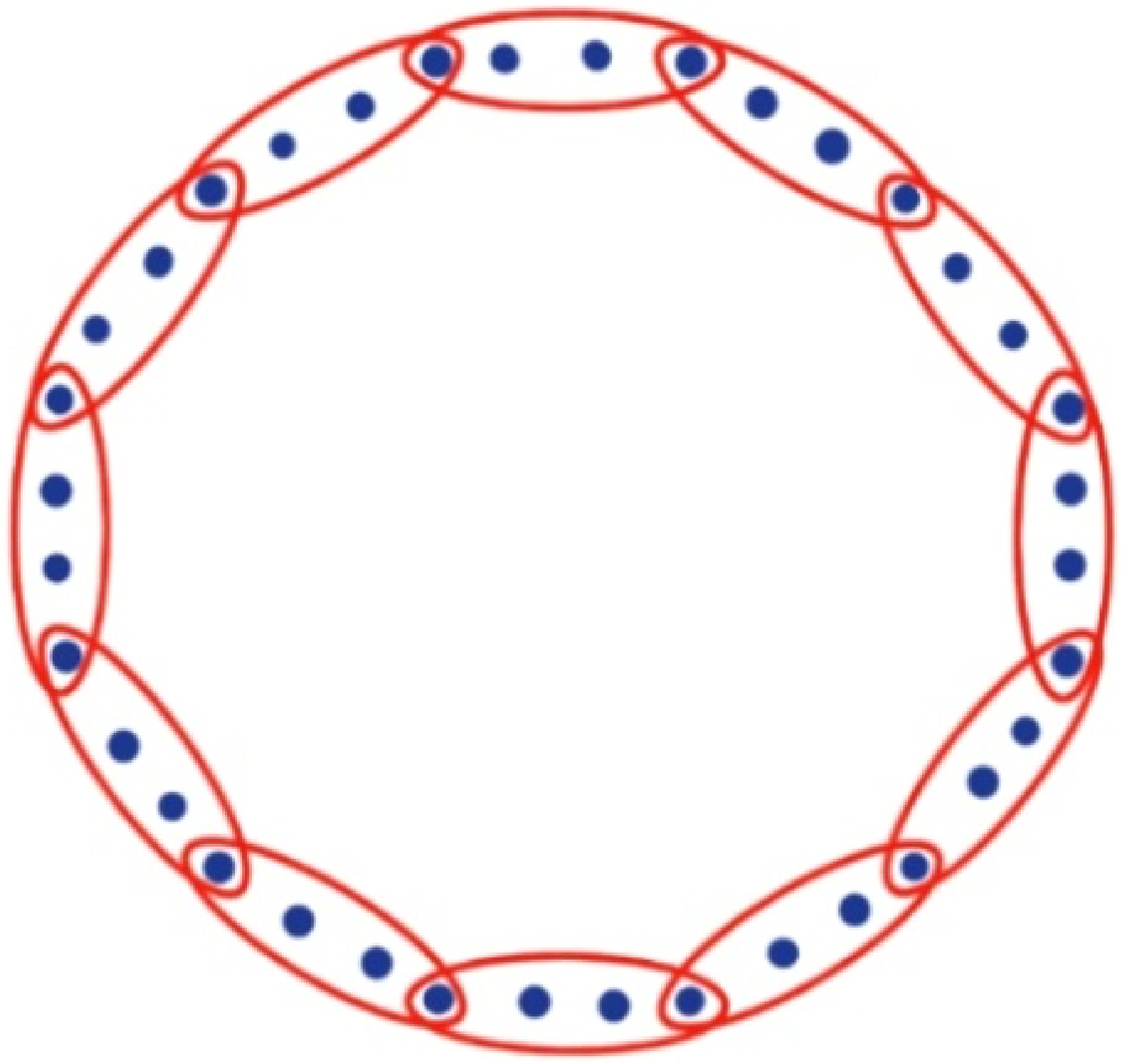}
\centerline{$(e) \ G^{4,12}_L $}
\end{minipage}%
\begin{minipage}[t]{0.25\linewidth}
\centering
\includegraphics[height=3cm]{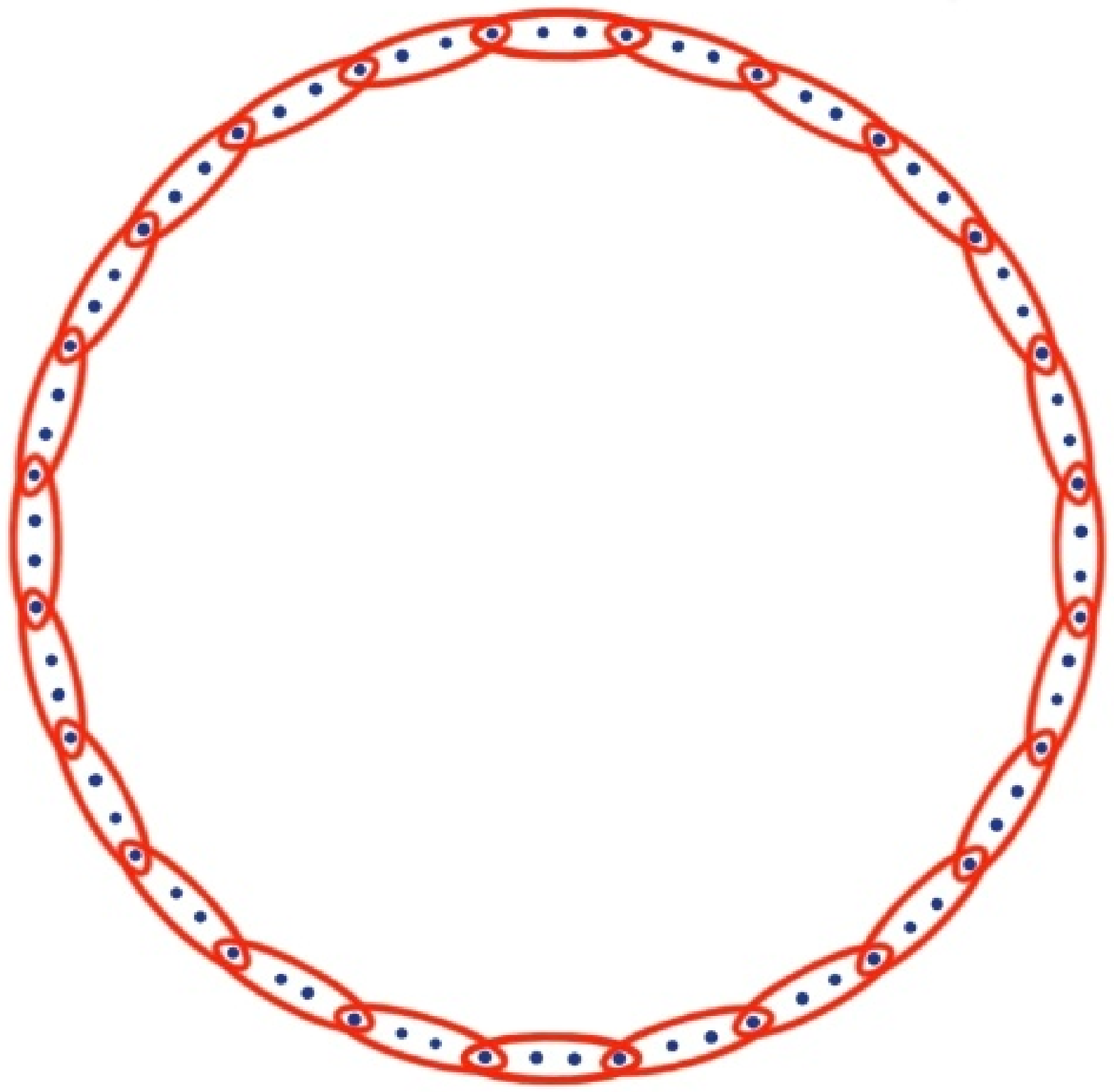}
\centerline{$(f) \ G^{4,24}_L $}
\end{minipage}
\caption{2-uniform loose cycles, the corresponding 4th power hypergraphs and other 4-uniform loose cycles.}
	\label{fig2}
\end{figure*}

$\mathbf{Example\ 3\ (Loose\ Cycle)}$. For a $k$-uniform hypergraph, if its vertex set $V=\{i_{1,1},\ldots,i_{1,k-1},i_{2,1},\ldots,i_{2,k-1},\ldots,i_{k,1},\ldots,i_{k,k-1}\}$ and its edge set
$E=\{(i_{1,1},\ldots,i_{1,k-1},i_{2,1}),(i_{2,1},\ldots,i_{2,k-1},i_{3,1}),\cdots,(i_{k,1},\ldots,i_{k,k-1},i_{1,1})\},$ then it is called a loose cycle. We denote the $k$-uniform loose cycles with $m$ edges as $G_L^{k,m}.$ Given a  graph $G$, if we add $k-2$ different vertices in its each edge, then we get its $k$th power hypergraph $G^k$. For example,  the 4-uniform loose cycles $G_L^{4,3}$ and $G_L^{4,6}$ in  Figure \ref{fig2} (c) and (d)  are the $4$th power hypergraph of the graphs  $G_L^{2,3}$ and $G_L^{2,6}$ in Figure \ref{fig2} (a) and (b) respectively. We list some conclusions about spectral radius  of  loose cycles and power hypergraphs below.
\begin{Proposition}[\cite{qi2017tensor}]\label{PowerHypergraph1}
If the spectral radius of the adjacency matrix of a graph $G$ is $\rho(A(G)),$ then the H-spectral radius of the adjacency tensor
  of its $k$th power hypergraph $\rho(\mathcal{A}(G^k))=\rho(A(G))^{\frac{2}{k}}$ .
\end{Proposition}
Since 2-uniform loose cycles $G_L^{2,m}$ are 2-regular, we obtain $\rho(A(G_L^{2,m}))=2$ from Proposition \ref{d-regular}. Therefore, based on Proposition \ref{PowerHypergraph1}, the  H-spectral radius of the adjacency tensor of the $4$th power hypergraph $G_L^{4,m}$ of $G_L^{2,m}$
\begin{equation}\label{HSpectralA}
 \rho(\mathcal{A}(G_L^{4,m}))=\rho(A(G_L^{2,m}))^{\frac{2}{4}}=\sqrt{2}, \ \text{for}\ \ m=3,4,5,\cdots.
\end{equation}

\begin{Proposition}[\cite{zhou2014some},\cite{hu2013cored}]\label{PowerHypergraph2}
  For a $d$-regular graph $G,$ the H-spectral radius of the signless Laplacian tensor of its $k$th power hypergraph $\rho(\mathcal{Q}(G^k))$ is the root of the equation
      \begin{equation}\label{rootofequation}
        (x-d)(x-1)^{\frac{k-2}{2}}-d=0.
      \end{equation}
If $k$ is even, then for any graph, the H-spectral radii of the signless Laplacian tensor and Laplacian tensor of its $k$th power hypergraph are equal. That is $\rho(\mathcal{Q}(G^k))=\rho(\mathcal{L}(G^k)).$
\end{Proposition}
When $k=4,d=2$ the root of equation \eqref{rootofequation} is 3. Thus from Proposition \ref{PowerHypergraph2},  the  H-spectral radius of the signless Laplacian tensor and Laplacian tensor of $G_L^{4,m}$
\begin{equation}\label{HSpectralQL}
 \rho(\mathcal{Q}(G_L^{4,m}))=\rho(\mathcal{L}(G_L^{4,m}))=3, \ \text{for} \ \ m=3,4,5, \cdots.
\end{equation}

%
%
%
%

We compute the largest H-eigenvalues of adjacency tensors $\mathcal{A}(G_L^{4,m})$ and Laplacian tensors $\mathcal{L}(G_L^{4,m})$ of the $4$-uniform loose cycles  in Figure \ref{fig2} (c), (d) and (e) when $m=3,6,12$. Table \ref{tab3} reports the results calculated by  PM, ACSA, HUOA and ACRCET. Obviously, the H-eigenvalues of $\mathcal{A}(G_L^{4,m}) $ and $\mathcal{L}(G_L^{4,m}) $ given by these four methods coincide with the theoretical results in \eqref{HSpectralA} and \eqref{HSpectralQL}. Compared with PM, ACSA and HUOA, ACRCET saves more than a half iterations and costs much less time. Because of the memory limitation of our laptop, PM, ACSA and HUOA are executable only under the condition $m \leq 48.$

We show the performance of ACRCET for computing the largest Z-eigenvalues of signless Laplacian tensors $\mathcal{Q}(G_L^{4,m})$ of 4-uniform loose cycles for different $m$ in Table \ref{tab4} . It can be seen that the ACRCET method is able to compute the  largest Z-eigenvalues of tensors with dimensions $n$ up to more than two thousands. Although the relationship between the  Z-spectral radius of a graph and the Z-spectral radius of its power hypergraph is not clear, it seems from the numerical results that a result about Z-spectral radius  similar to the conclusion in Proposition \ref{PowerHypergraph1} holds for $4$th power hypergraph of 2-uniform loose cycles.

\renewcommand{\arraystretch}{1.5} 
\begin{table}[tp]

  \centering
  \fontsize{10}{8}\selectfont
\scalebox{0.85}{
  \begin{threeparttable}
  \caption{Results for finding the largest H-eigenvalues of $\mathcal{A}(G_L^{4,m})$ and $\mathcal{L}(G_L^{4,m})$.}
  \label{tab3}

    \begin{tabular}{cc|ccc|ccc}
\hline
 & Algorithms&$\lambda_{\max}^H(\mathcal{A}(G_L^{4,m}))$&Iter'n&Time(s)&$\lambda_{\max}^H(\mathcal{L}(G_L^{4,m})$&Iter'n&Time(s)\cr
       \hline
  $m=3$ & PM&1.4142&4894&1.34&3.0000&4082&1.45\cr
   & ACSA&1.4142&2494&1.33&3.0000&3065&1.49\cr
   & HUOA&1.4142&2654&1.79&3.0000&2740&1.84\cr
    &ACRCET&{\bf 1.4142}&{\bf 532}&{\bf 0.47}&{\bf 3.0000}&{\bf 598}&{\bf 0.68}\cr
    \hline
   $m=6$ &PM&1.4142&14937&5.91&3.0000&13158&5.48\cr
   & ACSA&1.4142&8628&5.09&3.0000&8634&4.99\cr
    &HUOA&1.4142&4621&2.91&3.0000&5094&3.01\cr
   & ACRCET&{\bf 1.4142}&{\bf 808}&{\bf 0.92}&{\bf 3.0000}&{\bf 983}&{\bf 1.19}\cr
    \hline
   $m=12$ &PM&1.4142&48724&106.34&3.0000&46426&102.99\cr
   & ACSA&1.4142&18581&71.07&3.0000&16325&62.78\cr
    &HUOA&1.4142&8956&12.08&3.0000&10711&13.14\cr
   & ACRCET&{\bf 1.4142}&{\bf 1343}&{\bf 1.83}&{\bf 3.0000}&{\bf 1857}&{\bf 2.69}\cr
    \hline
    \end{tabular}

    \end{threeparttable}
}
\end{table}

\renewcommand{\arraystretch}{1.3} 
\begin{table}[tp]

  \centering
  \fontsize{10}{8}\selectfont
  \begin{threeparttable}
  \caption{Performance of ACRCET for finding the largest Z-eigenvalues of $\mathcal{Q}(G_L^{4,m})$.}
  \label{tab4}
    \begin{tabular}{cc|ccc}
\hline
  n&m&$\lambda_{\max}^Z(\mathcal{Q}(G_L^{4,m}))$&Iter'n&Time(s)\cr
    \hline
   9&3&2&350&0.42\cr
    18&6&2&340&0.69\cr
    36&12&2&635&1.10\cr
    72&24&2&586&1.76\cr
144&48&2&598&4.20\cr
288&96&2&690&15.06\cr
576&192&2&665&147.40\cr
1152&384&2&728&1020.07\cr
 2304&768&2&811&9745.41\cr
    \hline
    \end{tabular}
    \end{threeparttable}
\end{table}

\section{Conclusion}
In this paper we have used the adaptive cubic regularization method to compute extreme H- and Z-eigenvalues of even order symmetric tensors. We have established a fast computing skill, which has been proven effective in our test, for the matrix-valued products $\mathcal{T}\mathbf{x}^{r-2}$ of a vector $\mathbf{x}$ and a tensor $\mathcal{T}$ arising from a uniform hypergraph. Numerical experiments show that our ACRCET algorithm performs well for  even order symmetric tensor problems. Our next goal is to further study the cubic regularization method for  constrained optimization problems, and improve the computation efficiency of the cubic subproblem so that large scale problems can be calculated efficiently.

\bibliographystyle{plain}
\bibliography{ref}

\end{document}